\documentclass[12pt]{article}
\usepackage[margin=1in,top=1in,
            nohead,
            pdftex]{geometry}
\usepackage[T1]{fontenc}
\usepackage[utf8]{inputenc}
\usepackage{amsmath,amssymb,amsthm,mathtools,bm}
\usepackage{enumitem}
\usepackage{booktabs}
\usepackage{graphicx}
\usepackage{natbib}
\usepackage{bibunits}
\usepackage[colorlinks=true,citecolor=blue,urlcolor=blue,linkcolor=blue]{hyperref}

\usepackage{microtype}
\usepackage{xcolor}

\usepackage{algorithm}
\usepackage{algpseudocode}

\definecolor{darkgreen}{RGB}{0,100,0}

\newtheorem{theorem}{Theorem}
\newtheorem{lemma}{Lemma}
\newtheorem{proposition}{Proposition}
\newtheorem{assumption}{Assumption}

\theoremstyle{remark}
\newtheorem{remark}{Remark}
\newtheorem{example}{Example}

\allowdisplaybreaks
\setlist[itemize]{leftmargin=2.2em,itemsep=0.25em,topsep=0.25em}
\setlist[enumerate]{leftmargin=2.4em,itemsep=0.25em,topsep=0.25em}

\newcommand{\pr}{\mathbb P}
\newcommand{\E}{\mathbb E}

\newcommand{\one}{\mathbb{I}}

\newcommand{\floor}[1]{\left\lfloor #1\right\rfloor}

\title{Boot-AVT: Asymptotically Anytime-Valid Testing with Degenerate U-Statistics via Bootstrap}

\author{
	\begin{minipage}{0.96\textwidth}
		\centering
		\linespread{1}\selectfont
		Jixuan Liu\textsuperscript{1}, Long Feng\textsuperscript{1}, Liuhua Peng\textsuperscript{2}, and Zhaojun Wang\textsuperscript{1}\\[0.6em]
		{\small\textsuperscript{1}School of Statistics and Data Science, Nankai University,} 
		{\small Tianjin, China}\\[0.25em]
		{\small\textsuperscript{2}School of Mathematics and Statistics, University of Melbourne,} 
		{\small Victoria, Australia}
	\end{minipage}
}
\hypersetup{pdfauthor={Jixuan Liu, Long Feng, Liuhua Peng and Zhaojun Wang}}
\date{}

\begin{document}

\maketitle

\begin{abstract}

Many nonparametric tests are based on $U$-statistics that are degenerate under the null hypothesis. 
We develop Boot-AVT, a bootstrap procedure for asymptotically anytime-valid testing with these statistics. 
Using Brownian scaling and time inversion for the limiting Gaussian chaos process, we show that logarithmically penalized maxima over observed and future sample sizes have a common limiting null distribution. 
A quadratic multiplier bootstrap estimates the corresponding critical values. 
We also develop Boot-AVT-Multi for joint testing with multiple kernels. 
For both procedures, we establish uniform consistency of the bootstrap critical values, explicit bounds on size distortion, and asymptotic size $\alpha$ over an unbounded monitoring horizon as the monitoring start tends to infinity. 
Under fixed alternatives, we derive finite-horizon power bounds and establish power one with finite expected stopping times. 
Simulation studies evaluate finite-sample size and power in goodness-of-fit, two-sample and independence testing, and real-data analyses compare stopping times.

\end{abstract}

\section{Introduction}

Sequential hypothesis testing allows statistical evidence to be evaluated repeatedly as observations accumulate and permits early stopping when sufficient evidence against the null hypothesis is obtained. 
Applications include interim monitoring in clinical trials \citep{pocock1977}, online A/B experiments \citep{johari2022always}, and monitoring predictive performance after model deployment \citep{podkopaev2022tracking}. 
Sequential procedures can reduce the duration and cost of a study and need not require a prespecified final sample size. 
Unlike fixed-sample testing, sequential testing requires control of the overall type I error probability across monitoring times. 
Applying a fixed-sample test repeatedly at the nominal level does not, in general, provide this control \citep{johari2022always}.

Many nonparametric tests are based on $U$-statistics, which provide unbiased estimators of distributional functionals \citep{Hoeffding1948}. 
Examples include Cram\'er--von Mises goodness-of-fit statistics \citep{csorgo1996exact}, the maximum mean discrepancy (MMD) for two-sample testing \citep{gretton2012kernel}, and the Hilbert--Schmidt independence criterion (HSIC) for independence testing \citep{gretton2008kernel}. 
These statistics are degenerate under their respective null hypotheses and have non-Gaussian limit distributions. 
Permutation and bootstrap methods estimate their null distributions at a fixed sample size \citep{ArconesGine1992,gretton2008kernel,gretton2012kernel}. 
Fixed-sample distributional approximations do not by themselves establish sequential validity, which requires control of the probability that the statistic crosses a rejection boundary at any monitoring time.

Classical sequential procedures include sequential probability ratio tests \citep{wald1945sequential}, confidence sequences \citep{darling1967confidence}, and $\alpha$-spending \citep{landemets1983}. 
Nonnegative test martingales, supermartingales and $e$-processes provide nonasymptotic anytime-valid guarantees \citep{shafer2011test,howard2021timeuniform,ramdas2023game}. 
For kernel-based testing, \citet{balsubramani2016sequential} construct sequential tests using statistics formed from disjoint observation blocks, while \citet{manole2023martingale} derive concentration bounds for discrepancies and complete $U$-statistics using reverse-martingale methods. 
Betting procedures form wealth processes from predictable discrepancy witnesses \citep{shekhar2024nonparametric}, kernel dependence measures \citep{podkopaev2023skit}, or Stein discrepancies \citep{zhou2024sequential,martineztaboada2025sequential}. 
Ville's inequality then provides finite-sample anytime-valid error control without approximating the non-Gaussian null distribution of the original $U$-statistic.



Asymptotic approaches derive sequential boundaries through limiting stochastic processes. 
Invariance principles connect partial-sum boundary crossings to Brownian motion \citep{robbins1970boundary}, and recent work establishes time-uniform central limit theory and asymptotic confidence sequences for means and asymptotically linear estimators \citep{waudbysmith2024tuclt,GnettnerKirch2025}. 
For degenerate $U$-statistics, the null limit is governed by the second canonical term in the Hoeffding decomposition, and the limiting process is a Gaussian chaos whose distribution depends on an unknown spectrum \citep{DoblerKasprzakPeccati2022}. 
Thus, a first-order Gaussian approximation does not determine the relevant rejection boundary.

\citet{cai2026asymptotic} construct spectrally allocated Gaussian
chaos excursion (SAGE) boundaries for degree-two degenerate
$U$-statistics by allocating the significance level across
spectral coordinates and combining coordinatewise Gaussian
boundaries through a union bound. 
When the spectrum is unknown, the spectral quantities are estimated from the observations. 
Our approach directly calibrates the penalized maximum of the full limiting Gaussian chaos process. 
The main challenge is to estimate the required critical value from accumulated observations while maintaining validity when the estimate is updated after each observation.

We develop \textit{Boot-AVT} for asymptotically anytime-valid testing with degenerate $U$-statistics. 
Using Brownian scaling and time inversion, we show that the
logarithmically penalized maxima over observed and future sample
sizes share a common limiting null distribution. 
This relation allows the critical value for future monitoring
to be estimated from the accumulated observations. 
A quadratic multiplier bootstrap based on the empirical second
canonical function estimates its null distribution without explicit eigenvalue estimation. 
Within each replication, the same multiplier sequence is used across sample sizes to account for the dependence among the statistics entering the maximum. 
We establish uniform consistency of the bootstrap critical values
over all subsequent sample sizes. 
The resulting test has asymptotic size $\alpha$ as the monitoring
start tends to infinity, with the critical value re-estimated
after each observation. 


The main contributions are as follows.
\begin{enumerate}[label=(\roman*)]
    \item We develop Boot-AVT for sequential tests based on degenerate $U$-statistics of arbitrary fixed degree. 
    Using Brownian scaling and time inversion of the limiting Gaussian chaos process, we establish a common limiting null distribution for the logarithmically penalized maxima over observed and future sample sizes. 
    This relation allows the critical value for monitoring over an unbounded horizon to be estimated from the observations available at each sample size.

    \item We construct a quadratic multiplier bootstrap based on the empirical second canonical function, without explicit eigenvalue estimation. 
    We prove almost-sure uniform consistency of the estimated critical values over all subsequent sample sizes and establish that the overall type I error probability converges to $\alpha$ as the monitoring start tends to infinity, with the critical value re-estimated after each observation. 
    We derive explicit bounds on size distortion that separate statistical approximation and Monte Carlo errors.
    Under fixed alternatives, we also obtain finite-horizon power bounds and establish power one with finite expected stopping times. 

    \item We develop Boot-AVT-Multi for joint monitoring of a fixed finite collection of kernels or feature representations. 
    Common multipliers account for dependence across candidate statistics and sample sizes, permitting direct calibration of the maximum over both kernels and monitoring times. 
    We establish asymptotically exact time-uniform error control for the joint test and extend the fixed-alternative power guarantees to settings in which at least one candidate functional exceeds its null value.
\end{enumerate}

Simulations evaluate finite-sample size and power in goodness-of-fit, two-sample and independence testing. 
Applications to audiometry and breast cancer data compare the stopping times of the proposed and competing tests. 

The remainder of the paper is organized as follows. Section~\ref{sec:preliminary} introduces the sequential error criterion, degenerate $U$-statistics and the motivating kernel statistics. 
Section~\ref{sec:methodology} presents the boundary construction based on time inversion and a quadratic multiplier bootstrap, establishes the theoretical guarantees, and describes the online implementation. The extensions to multiple kernels are proposed in Section~\ref{sec:joint-finite-family}. Sections~\ref{sec:simulations} and~\ref{sec:realdata} report simulations and real-data applications, respectively. 


\section{Preliminaries}\label{sec:preliminary}

Let $X_1,X_2,\ldots$ be i.i.d.~observations from an unknown distribution $P$ on $\mathcal X\subseteq\mathbb R^d$, and let $\theta=T(P)\in\mathbb R$ be a functional of interest. 
For a specified null value $\theta_0$, we consider
\begin{equation}
    \notag
    H_0:\theta=\theta_0 \qquad\text{versus}\qquad H_1:\theta>\theta_0.
\end{equation}

We assume that $\theta=\E\{h(X_1,\ldots,X_q)\}$, where $h:\mathcal X^q\to\mathbb R$ is a symmetric kernel of fixed degree $q\geq2$ satisfying $\E|h(X_1,\ldots,X_q)|<\infty$.
An unbiased estimator of $\theta$ based on the first $n$ observations is the $U$-statistic
\begin{equation}
    \notag 
    U_{n,q}(h) = \binom nq^{-1}\sum_{1\le i_1<\cdots<i_q\le n}h(X_{i_1},\ldots,X_{i_q})
\end{equation}
for $n\geq q$.
The first two canonical functions in the Hoeffding decomposition \citep{Hoeffding1948,Serfling1980} of the $U$-statistic are
$h^{(1)}(x)=\E\{h(x,X_2,\ldots,X_q)\}-\theta$ and
\begin{equation}
    \notag 
    h^{(2)}(x,y) = \E\{h(x,y,X_3,\ldots,X_q)\}-\theta-h^{(1)}(x)-h^{(1)}(y).
\end{equation}
For $q=2$, the expectation in the second formula is understood as $h(x,y)$. 
We assume that, under $H_0$, $h^{(1)}(X_1)=0$ almost surely and
$0<\E[\{h^{(2)}(X_1,X_2)\}^2]<\infty$.
Thus, the $U$-statistic is degenerate of order two \citep{Serfling1980}. 
Let $a_q=q(q-1)/2$, and define the normalized statistic
\begin{equation}
    \notag 
    T_n(\theta) = n\{U_{n,q}(h)-\theta\}/a_q.
\end{equation}
In the following examples, $\theta$ is a nonnegative discrepancy that equals zero precisely under the null hypothesis, so the alternative is equivalent to $\theta>0$.

\begin{example}[Goodness-of-fit test]\label{ex:prelim-gof}
    Let $F$ be the distribution function of the univariate observations $X_i$, and let $F_0$ be a specified continuous distribution function.
    Consider
    \begin{equation}
        \notag 
        H_0:F=F_0 \qquad\text{versus}\qquad H_1:F\ne F_0.
    \end{equation}
    The Cram\'er--von Mises discrepancy is 
    $\theta=\int_{\mathbb R}\{F(t)-F_0(t)\}^2\,dF_0(t)$,
    which is zero if and only if $F=F_0$  \citep{csorgo1996exact}.
    It admits the representation $\theta=\E\{h(X_1,X_2)\}$ with the symmetric degree-two kernel
    $h(x,y) = \int_{\mathbb R}\{\one(x\le t)-F_0(t)\}\{\one(y\le t)-F_0(t)\}\,dF_0(t)$,
    where $\one(\cdot)$ denotes the indicator function.
    Under $H_0$, the kernel is degenerate of order two.
\end{example}

\begin{example}[Two-sample MMD test]\label{ex:prelim-mmd}
    Let $X_i=(Y_i,Z_i)$, $i\geq1$, be i.i.d.\ pairs with marginal distributions $P_Y$ and $P_Z$. 
    Consider
    \begin{equation}
        \notag 
        H_0:P_Y=P_Z \qquad\text{versus}\qquad H_1:P_Y\ne P_Z.
    \end{equation}
    For a kernel $K$, the squared maximum mean discrepancy \citep[MMD]{gretton2012kernel} is
    $\theta=\E\{K(Y,Y')\}+\E\{K(Z,Z')\}-2\E\{K(Y,Z')\}$,
    where $X=(Y,Z)$ and $X'=(Y',Z')$ are independent copies of $X_1$.
    If $K$ is characteristic, $\theta=0$ if and only if $P_Y=P_Z$.
    An unbiased estimator based on $X_1,\ldots,X_n$ is $U_{n,2}(h)$, with
    $h(x,x') = K(y,y')+K(z,z')-K(y,z')-K(y',z)$,
    where $x=(y,z)$ and $x'=(y',z')$.
    Under $H_0$, $h^{(1)}=0$ and $h^{(2)}=h$, so the $U$-statistic is degenerate of order two provided $\E[\{h(X_1,X_2)\}^2]>0$.
\end{example}

\begin{example}[HSIC Independence test]\label{ex:prelim-hsic}
    Let $X_i=(Y_i,Z_i)$, $i\ge1$, be i.i.d.~paired observations with joint distribution $P_{YZ}$ and marginals $P_Y$ and $P_Z$.
    Consider
    \begin{equation}
        \notag 
        H_0:P_{YZ}=P_Y\otimes P_Z \qquad\text{versus}\qquad H_1:P_{YZ}\ne P_Y\otimes P_Z,
    \end{equation}
    where $\otimes$ denotes a product measure.
    For two prespecified bounded, continuous, characteristic kernels $K_1$ and $K_2$, the Hilbert--Schmidt independence criterion \citep[HSIC]{gretton2008kernel} is 
    $\theta=\E\{K_1(Y,Y')K_2(Z,Z')\}+\E\{K_1(Y,Y')\}\E\{K_2(Z,Z')\}-2\E\{K_1(Y,Y')K_2(Z,Z'')\}$,
    where $X=(Y,Z)$, $X'=(Y',Z')$ and $X''=(Y'',Z'')$ are independent copies of $X_1$.
    This criterion vanishes if and only if $Y$ and $Z$ are independent.
    An unbiased estimator of $\theta$ is the degree-four $U$-statistic with symmetric kernel 
    \begin{equation}
        \notag 
        h(x_1,x_2,x_3,x_4) = \frac{1}{24}\sum_{\upsilon\in\mathcal S_4}K_1(y_{\upsilon_1},y_{\upsilon_2})\bigl\{K_2(z_{\upsilon_1},z_{\upsilon_2}) +K_2(z_{\upsilon_3},z_{\upsilon_4})-2K_2(z_{\upsilon_1},z_{\upsilon_3})\bigr\},
    \end{equation}
    where $x_i=(y_i,z_i)$, $\upsilon=(\upsilon_1,\upsilon_2,\upsilon_3,\upsilon_4)$, and $\mathcal S_4$ is the set of permutations of $\{1,2,3,4\}$.
    Under $H_0$, $h^{(1)}=0$, and the kernel is degenerate of order two provided $\E[\{h^{(2)}(X_1,X_2)\}^2]>0$.
\end{example}

For sequential testing, fix a nominal level $\alpha\in(0,1)$,
an initial sample size $m\geq q$ and a parameter $c>0$. 
Monitoring begins at sample size $m$ and continues after each additional observation. 
At each sample size $n\geq m$, let $\kappa_{1-\alpha,c}(n)$ be a
critical value computed from $X_1,\ldots,X_n$. 
We consider the rejection rule and corresponding stopping time
\begin{equation}
    \notag 
    \phi_{m,n} = \one\{T_n(\theta_0)-c\log(n/m)>\kappa_{1-\alpha,c}(n)\} \qquad \text{and} \qquad 
    \tau_m=\inf\{n\ge m:\phi_{m,n}=1\},
\end{equation}
with $\inf\varnothing=\infty$.
Let $\pr_0$ denote probability under the null hypothesis.
The time-uniform (or overall) type I error probability is
\begin{equation}
    \notag 
    \pr_0(\tau_m<\infty) = \pr_0\left(\exists~n\ge m:\phi_{m,n}=1\right),
\end{equation}
which accounts for rejection at any time over the entire monitoring horizon.
The procedure with a sequence of tests $\{\phi_{m,n}\}_{n\geq m}$ is asymptotically anytime valid at level $\alpha$ if
\begin{equation}
    \notag
    \limsup_{m\to\infty}\pr_0(\tau_m<\infty)=\limsup_{m\to\infty}\pr_0\left(\exists~n\geq m:\phi_{m,n}=1\right)\le\alpha.
\end{equation}
We seek asymptotically exact error control, or asymptotic size $\alpha$, in the sense that
\begin{equation}
    \notag 
    \lim_{m\to\infty}\pr_0(\tau_m<\infty)=\lim_{m\to\infty}\pr_0\left(\exists~n\geq m:\phi_{m,n}=1\right)=\alpha.
\end{equation}
In addition, we derive explicit bounds on $|\pr_0(\tau_m<\infty)-\alpha|$.
A bound tending to zero establishes asymptotic size $\alpha$,
while its rate quantifies how rapidly the overall type I error
probability converges to $\alpha$ as $m\to\infty$.


\section{Methodology}\label{sec:methodology}

Boot-AVT approximates a logarithmic rejection boundary using the null distribution of a penalized $U$-statistic maximum. 
We first study the limiting process and its boundary-crossing probability, then construct the multiplier bootstrap and establish its sequential guarantees.

For the statistic $T_n(\theta)$ defined in Section~\ref{sec:preliminary}, consider the rescaled $U$-process
\begin{equation}\label{eq:rescaled-U-process}
    T_{\lfloor nt\rfloor}(\theta) = \lfloor nt\rfloor\{U_{\lfloor nt\rfloor,q}(h)-\theta\}/a_q, \qquad t\geq q/n.
\end{equation}
Let $\{\lambda_\ell\}_{\ell\ge 1} \subset\mathbb R$ denote the eigenvalues of the integral operator with kernel $h^{(2)}(\cdot,\cdot)$, and let $\{\psi_\ell(\cdot)\}_{\ell\ge 1}$ be the corresponding orthonormal eigenfunctions in $L^2(P)$.
Then
$h^{(2)}(x,y)=\sum_\ell\lambda_\ell\psi_\ell(x)\psi_\ell(y)$ 
in $L^2(P\times P)$, where $\E\{\psi_\ell(X_1)\}=0$ and
$\E\{\psi_\ell(X_1)\psi_r(X_1)\}=\one(\ell=r)$ for $\ell,r\in\mathbb N=\{1,2,\ldots\}$.
We impose the following spectral and moment conditions.

\begin{assumption}\label{ass:front-spectrum}
    The eigenvalues satisfy $0<\Lambda_1=\sum_{\ell\ge 1}|\lambda_\ell|<\infty$.
\end{assumption}

\begin{assumption}\label{ass:front-kernel}
    There exists a constant $p>2$ such that $\E|h(X_1,\ldots,X_q)|^p<\infty$.
\end{assumption}

\begin{remark}\label{rem:front-spectrum-examples}
    Assumption~\ref{ass:front-spectrum} also appears in \citet{cai2026asymptotic,SejdinovicEtAl2013}. 
    The condition holds for the Cram\'er--von Mises kernel and for MMD and HSIC based on Gaussian, Laplace or Mat\'ern kernels under their respective null hypotheses \citep{gretton2012kernel,SongEtAl2012,SejdinovicEtAl2013}.
\end{remark}

Under $H_0$, the Hoeffding decomposition gives $T_{\lfloor nt\rfloor}(\theta)=\lfloor nt\rfloor U_{\lfloor nt\rfloor,2}(h^{(2)})+o_p(1)$, and the quadratic term $\lfloor nt\rfloor U_{\lfloor nt\rfloor,2}(h^{(2)})$ determines the weak limit of $T_{\lfloor nt\rfloor}(\theta)$.
Let $W_\ell$ be independent standard Brownian motions and define the Gaussian chaos process \citep{DoblerKasprzakPeccati2022}
\begin{equation}
    \notag
    Q(t)=\sum_{\ell\ge 1}\lambda_\ell\left\{W_\ell^2(t)/t-1\right\},\qquad t>0.
\end{equation}
Let $D[\delta_L,\delta_U]$ denote the space of right-continuous functions with left limits on $[\delta_L,\delta_U]$.

\begin{lemma}\label{lem:1}
    Under $H_0$, suppose Assumptions~\ref{ass:front-spectrum}--\ref{ass:front-kernel} hold. 
    Then, for every fixed $0<\delta_L<\delta_U<\infty$,
    \begin{equation*}\label{eq:front-U-FCLT}
        \{T_{\floor{nt}}(\theta):\delta_L\leq t\leq \delta_U\}\Rightarrow\{Q(t):\delta_L\leq t\leq \delta_U\}
    \end{equation*}
    in $D[\delta_L,\delta_U]$ as $n\to\infty$.
\end{lemma}

Lemma~\ref{lem:1} gives weak convergence of the rescaled degenerate $U$-process to the Gaussian chaos process on each fixed time interval $[\delta_L,\delta_U]$. 
Sequential testing, however, requires control of the overall
probability of rejecting the null hypothesis when the test
is applied at every sample size $n\geq m$. 
We therefore need to approximate the null distribution of a supremum over an unbounded time interval and estimate the corresponding critical value from the observed sample. 
Boot-AVT determines this critical value from the distribution of a logarithmically penalized maximum of $Q(t)$.

\subsection{Boundary construction by time inversion}\label{sec:oracle-Gaussian-critical-value}

We first construct an oracle boundary for $\{Q(n)\}_{n\geq m}$,
treating the eigenvalues $\{\lambda_\ell\}_{\ell\ge 1}$ as known.
For fixed tuning parameters, Gaussian-mixture boundaries for
sample means have order $\sqrt{\log n/n}$, so their squares have
order $\log n/n$
\citep{robbins1970boundary,howard2021timeuniform}.
We consider the same order for the boundary of
$U_{n,q}(h)-\theta_0$.
On the scale of $T_n(\theta_0)$ and its limiting process,
this corresponds to logarithmic growth.

For fixed $c>0$ and an integer $m\geq1$, consider the boundary
$\kappa+c\log(n/m)$ for $Q(n)$ over $n\geq m$. 
To determine $\kappa$, we first consider the corresponding
boundary for the continuous-time process $Q(t)$ over $t\ge m$  and, for $\alpha\in(0,1)$, require
$\pr\bigl(\sup_{t\geq m}\{Q(t)-c\log(t/m)\}>\kappa\bigr)=\alpha$.
The resulting boundary is then applied at integer monitoring
times.

To estimate this critical value without explicitly estimating the
eigenvalues, we seek a maximum based on the first $n$ observations whose null distribution approximates that of the continuous-time supremum.
Brownian scaling and time inversion relate this supremum over
$[m,\infty)$ to one over $(0,n]$.
Under Assumption~\ref{ass:front-spectrum}, for every $m>0$ and integer $n\geq m$,
\begin{align}\label{eq:front-time-inversion}
    \notag 
    & \sup_{t\geq m}\{Q(t)-c\log(t/m)\}
    \overset d=\sup_{u\geq1}\{Q(u)-c\log u\} \\
    \overset d= & \sup_{0<s\leq1}\{Q(s)+c\log s\} 
    \overset d= \sup_{0<t\leq n}\{Q(t)+c\log(t/n)\}.
\end{align}
The first and third equalities follow from Brownian scaling. 
The second follows from the time-inversion property that
$\{sW_\ell(1/s):s>0\}$ and $\{W_\ell(s):s>0\}$ have the same
process distribution \citep{RevuzYor1999}. 
The last supremum involves only times up to $n$. 
Its $U$-statistic counterpart is a maximum of $T_k(\theta_0)+c\log(k/n)$ over $q\leq k\leq n$, which can be
computed from the first $n$ observations. 
This relation determines the form of the maximum used to estimate the critical value for future monitoring.


Let $\Phi_c(s)=Q(s)+c\log s$ and define the penalized maximum
$H_c=\sup_{0<s\leq1}\Phi_c(s)$. 
The logarithmic term controls the behavior near zero. 
Under Assumption~\ref{ass:front-spectrum}, $Q(s)/|\log s|\to0$ almost surely as $s\downarrow0$, so $\Phi_c(s)\to-\infty$ almost surely. 
Together with the continuity of $Q(s)$ on $(0,1]$, this implies that $H_c$ is finite and the supremum is attained almost surely. 
In addition, its distribution function $F_c(x)=\pr(H_c\leq x)$ is continuous
and strictly increasing on the interior of its support. 
Consequently, its $(1-\alpha)$-quantile,
\begin{equation*}
    \kappa_{1-\alpha,c} = \inf\{x\in\mathbb R:F_c(x)\geq1-\alpha\},
\end{equation*}
satisfies $F_c(\kappa_{1-\alpha,c})=1-\alpha$.
Equation~\eqref{eq:front-time-inversion} therefore yields
\begin{align}\label{eq:exact-continuous-coverage}
    \notag 
    & \pr\left(\sup_{t\geq m}\{Q(t)-c\log(t/m)\} \leq \kappa_{1-\alpha,c}\right) \\
    = & \pr\left(Q(t)-c\log(t/m)\leq\kappa_{1-\alpha,c}
    \text{ for all }t\geq m\right)=1-\alpha.
\end{align}
Thus, a single critical value $\kappa_{1-\alpha,c}$,
independent of $m$, determines the rejection boundary over the entire
continuous-time monitoring interval $[m,\infty)$.
For integer monitoring times, the following theorem establishes
simultaneous coverage of at least $1-\alpha$ and gives a rate of
convergence of the coverage probability to $1-\alpha$ as
$m\to\infty$.

\begin{theorem}\label{thm:oracle-Gaussian chaos-boundary}
    Suppose Assumption~\ref{ass:front-spectrum} holds.
    For fixed $\alpha\in(0,1)$ and $c>0$, there exists a finite
    constant $C>0$ such that, for every integer $m\geq1$,
    \begin{align}
        \notag 
        0 \leq \pr\left(Q(n)-c\log(n/m)\leq\kappa_{1-\alpha,c} \text{ for all } n\geq m\right)-(1-\alpha) \leq C\{\log(m+1)/\sqrt m\}^{1/5}.
    \end{align}
\end{theorem}

The convergence rate of the coverage probability of the rejection boundary to $1-\alpha$ is $O(\{\log(m+1)/\sqrt m\}^{1/5})$. 
This error bound reflects discretization of the monitoring interval. 
It combines a bound on the increments of $Q$ between consecutive
integers with the local $1/5$-H\"older bound for $F_c$.


To approximate the last supremum in
\eqref{eq:front-time-inversion} using the integer times
$1,\ldots,n$, define
\begin{equation}
    \notag 
    {G}_{c,n}=\max_{1\leq k\leq n}\{Q(k)+c\log(k/n)\},
\end{equation}
and let $\kappa_{1-\alpha,c}^{\mathrm{grid}}(n)$ be its $(1-\alpha)$-quantile. 
It can be shown that $G_{c,n}$ converges to $H_c$ in distribution and $\kappa_{1-\alpha,c}^{\mathrm{grid}}(n)$ converges to $\kappa_{1-\alpha,c}$ as $n\to\infty$.
The next theorem shows that replacing the continuous-time critical value by its finite-grid counterpart still preserves the asymptotic time-uniform error control, with a polynomial convergence rate in the monitoring start $m$.

\begin{theorem}\label{thm:finite-grid-Gaussian chaos-boundary}
    Suppose Assumption~\ref{ass:front-spectrum} holds. 
    Then, for fixed $\alpha\in(0,1)$ and $c>0$, as $m\to\infty$,
    \begin{equation*}
        \pr\left(Q(n)-c\log(n/m)\leq\kappa_{1-\alpha,c}^{\mathrm{grid}}(n) \text{ for all } n\geq m
        \right) = 1-\alpha+O(m^{-\min\{1/40,c/(32e\Lambda_1)\}}).
    \end{equation*}
\end{theorem}

We next extend the boundary construction to the $U$-statistic
process. For $m\geq q$, define
\begin{equation*}
    J_{c,m}^{U} = \sup_{k\geq m}\{T_k(\theta)-c\log(k/m)\}.
\end{equation*}
Under $H_0$, the overall rejection probability for the boundary
$\kappa+c\log(k/m)$ is $\pr_0(J_{c,m}^{U}>\kappa)$.
To estimate the critical value from the observed sample, the
time-inversion identity in \eqref{eq:front-time-inversion}
suggests the maximum
\begin{equation*}
    H_{c,n}^{U} = \max_{q\leq k\leq n}\{T_k(\theta)+c\log(k/n)\},
    \qquad n\geq q,
\end{equation*}
which is computed from $X_1,\ldots,X_n$ under $H_0$.
The following theorem establishes that $H_{c,n}^{U}$ and
$J_{c,m}^{U}$ have the same limiting null distribution.
This common limit allows the critical value for future monitoring
to be estimated by approximating the null distribution of
$H_{c,n}^{U}$, as developed in
Section~\ref{sec:bootstrap-critical-value}.

\begin{theorem}\label{thm:oracle-Gaussian chaos-U-boundary}
Under $H_0$, suppose Assumptions~\ref{ass:front-spectrum} and~\ref{ass:front-kernel} hold. Let $\eta=48^{-1}\min\{(1-2/p)^2,1/9\}$, $\alpha\in(0,1)$ and $c>0$. Then $H_{c,n}^{U}$ converges to $H_c$ in distribution as $n\to\infty$ and $J_{c,m}^{U}$ converges to $H_c$ in distribution as $m\to\infty$. Moreover, as $m\rightarrow\infty$,
\begin{equation*}
\pr_0\left(
T_n(\theta_0)-c\log(n/m)\leq
\kappa_{1-\alpha,c}^{\mathrm{grid}}(n)
\text{ for all }n\geq m
\right)
=1-\alpha+O\left(m^{-\min\{\eta,c/(32e\Lambda_1)\}}\right).
\end{equation*}
\end{theorem}

Theorem~\ref{thm:oracle-Gaussian chaos-U-boundary} implies that rejecting $H_0$ when $U_{n,q}(h)-\theta_0$ exceeds $a_qn^{-1}\{\kappa_{1-\alpha,c}^{\mathrm{grid}}(n)+c\log(n/m)\}$
gives a test with asymptotic size $\alpha$.
For fixed $m$ and $c>0$, the boundary has order $\log n/n$, matching the order of a squared Gaussian-mixture boundary for a mean \citep{robbins1970boundary,howard2021timeuniform}.
The coverage-error bound combines the finite-grid approximation error with the error in approximating the $U$-statistic maximum over the unbounded monitoring horizon. 
Both terms have polynomial rates in $m$.

The finite-grid critical values
$\kappa_{1-\alpha,c}^{\mathrm{grid}}(n)$ and their limit
$\kappa_{1-\alpha,c}$ still depend on the unknown spectrum of the
second canonical function. 
Rather than estimating the eigenvalues, we directly approximate
the null distribution of $H_{c,n}^{U}$. 
The time-inversion identity in \eqref{eq:front-time-inversion} and the common limit of $J_{c,m}^{U}$ and
$H_{c,n}^{U}$ in Theorem~\ref{thm:oracle-Gaussian chaos-U-boundary} justify this choice: the $(1-\alpha)$-quantile of
$H_{c,n}^{U}$ converges to $\kappa_{1-\alpha,c}$. 
This relation allows the critical value for future monitoring
to be estimated by bootstrapping the maximum over observed sample
sizes. 
Section~\ref{sec:bootstrap-critical-value} develops a quadratic
multiplier bootstrap for this purpose.

\subsection{Multiplier bootstrap and asymptotic validity}\label{sec:bootstrap-critical-value}

Under $H_0$, the leading Hoeffding term of $T_k(\theta_0)$ is
$kU_{k,2}(h^{(2)})$. 
Boot-AVT uses a quadratic multiplier bootstrap based on an estimate of $h^{(2)}(\cdot,\cdot)$ to approximate the null distribution of $H_{c,n}^{U}$ and estimate $\kappa_{1-\alpha,c}$ \citep{Han2022}.
The multipliers satisfy the following condition.

\begin{assumption}\label{ass:front-multiplier}
    The multipliers $\xi_1,\xi_2,\ldots$ are i.i.d., independent of the observations, and satisfy $\E(\xi_1)=0$, $\E(\xi_1^2)=1$ and $\E|\xi_1|^p<\infty$.
\end{assumption}

We estimate $h^{(2)}(\cdot,\cdot)$ by replacing its defining expectations with averages over the remaining observations.
For $1\leq i<j\leq k$, define
\begin{align}\label{eq:empirical-conditional-kernel}
    \notag 
    \widehat g_{k}^{(2)}(X_i,X_j)
    &=\binom{k-2}{q-2}^{-1}
    \sum_{\substack{1\leq r_1<\cdots<r_{q-2}\leq k\\
                     r_\ell\notin\{i,j\},\ 1\leq\ell\leq q-2}}
    h(X_i,X_j,X_{r_1},\ldots,X_{r_{q-2}}),\\
    \widehat g_{k}^{(1)}(X_i)
    &=\binom{k-1}{q-1}^{-1}
    \sum_{\substack{1\leq r_1<\cdots<r_{q-1}\leq k\\r_\ell\ne i,\ 1\leq\ell\leq q-1}}
    h(X_i,X_{r_1},\ldots,X_{r_{q-1}}).
\end{align}
For $q=2$, $\widehat g_{k}^{(2)}(X_i,X_j)=h(X_i,X_j)$.
The estimator of the second canonical function is
\begin{align}\label{eq:empirical-second-canonical}
    \widehat h_k^{(2)}(X_i,X_j) = \widehat g_{k}^{(2)}(X_i,X_j)
    -\widehat g_{k}^{(1)}(X_i)-\widehat g_{k}^{(1)}(X_j)+U_{k,q}(h).
\end{align}

For each bootstrap replication $b$, generate an independent copy
$\xi_1^{(b)},\xi_2^{(b)},\ldots$ of the multiplier sequence and define
\begin{equation}
    \notag 
    \widetilde T_{k,2}^{*,b}
    =\frac{2}{k-1}\sum_{1\leq i<j\leq k}
    \xi_i^{(b)}\xi_j^{(b)}\widehat h_k^{(2)}(X_i,X_j),
    \qquad k\geq q.
\end{equation}
The bootstrap counterpart of $H_{c,n}^{U}$ is
\begin{equation}\label{eq:bootstrap-maximum}
    H_{c,n}^{*,b} = \max_{q\leq k\leq n}\left\{\widetilde T_{k,2}^{*,b}+c\log(k/n)\right\}.
\end{equation}
Within each replication, the same multiplier sequence is used across sample sizes to approximate the joint null distribution of the statistics entering the maximum. 
Each $\widetilde T_{k,2}^{*,b}$ uses $\widehat h_k^{(2)}(\cdot,\cdot)$ computed
from $X_1,\ldots,X_k$ and is retained as $n$ increases.
The parameter $c>0$ is fixed throughout monitoring; its choice is discussed in Remark~\ref{rem:choice-c} in Section~\ref{sec:online-implementation}.

Let
$F_n^*(x)=\pr(H_{c,n}^{*,1}\leq x\mid X_1,\ldots,X_n)$
denote the conditional distribution function of $H_{c,n}^{*,1}$.
For a prespecified deterministic sequence of positive integers
$\{B_n\}_{n\geq q}$, let
$\widehat\kappa_{1-\alpha,c}^{*,(B_n)}(n)$ be the empirical quantile of $H_{c,n}^{*,1},\ldots,H_{c,n}^{*,B_n}$. Equivalently,
\begin{equation}
    \notag 
    \widehat\kappa_{1-\alpha,c}^{*,(B_n)}(n) = \inf\{x\in\mathbb R:\widehat F_{n,B_n}^*(x)\geq1-\alpha\},
\end{equation}
where 
$\widehat F_{n,B_n}^*(x)=B_n^{-1}\sum_{b=1}^{B_n}\one(H_{c,n}^{*,b}\leq x)$ which approximates $F_n^*(x)$. 
We use $\widehat\kappa_{1-\alpha,c}^{*,(B_n)}(n)$ to define the Boot-AVT rejection rule as
\begin{equation}
    \notag
    \phi_{m,n}^{(B_n)}(\theta_0) = \one\left\{T_n(\theta_0)-c\log(n/m)>\widehat\kappa_{1-\alpha,c}^{*,(B_n)}(n)\right\},
\end{equation}
and the corresponding stopping time as
$\tau_m^{(B_n)}=\inf\{n\geq m:\phi_{m,n}^{(B_n)}(\theta_0)=1\}$. 
The following theorem establishes consistency of the bootstrap approximation and the estimated critical values, and gives an explicit bound on the overall type I error probability. 

\begin{theorem}\label{thm:front-final-bootstrap}
Under $H_0$, suppose Assumptions~\ref{ass:front-spectrum}--\ref{ass:front-multiplier} hold. Fix $\alpha\in(0,1)$ and $c>0$, and suppose that $B_n/\log(n+1)\to\infty$. Let $\underline b_m=\inf_{n\geq m}\{B_n/\log(n+1)\}$. Then $\sup_{n\geq m}|\widehat\kappa_{1-\alpha,c}^{*,(B_n)}(n)-\kappa_{1-\alpha,c}|\to0$ almost surely as $m\to\infty$. Moreover, there exists a finite constant $C>0$ such that, for all sufficiently large integers $m$,
\begin{align}
&\left|\pr_0\left(
T_n(\theta_0)-c\log(n/m)\leq
\widehat\kappa_{1-\alpha,c}^{*,(B_n)}(n)
\text{ for all }n\geq m\right)-(1-\alpha)\right|\notag\\
\leq&
C\{\log\log(m+1)\}^{1/2}\{\log(m+1)\}^{-p/2}
+(2/\underline b_m)^{1/2}+2/(3m^3).
\label{eq:bootstrap-size-error}
\end{align}
Consequently,
\begin{equation*}
\begin{aligned}
&\pr_0\left(\exists~n\geq m:\phi_{m,n}^{(B_n)}(\theta_0)=1\right)
=\pr_0(\tau_m^{(B_n)}<\infty)\\
=&\alpha+O\left(
\{\log\log(m+1)\}^{1/2}\{\log(m+1)\}^{-p/2}
+\underline b_m^{-1/2}\right).
\end{aligned}
\end{equation*}
\end{theorem}
The uniform convergence of the estimated critical values, together with the limit of $J_{c,m}^{U}$ in Theorem~\ref{thm:oracle-Gaussian chaos-U-boundary}, yields asymptotic size $\alpha$ when the critical value is re-estimated
at every sample size. 
The first term in \eqref{eq:bootstrap-size-error} bounds the error from the process and bootstrap approximations, including estimation of $h^{(2)}(\cdot,\cdot)$.
The remaining terms bound the Monte Carlo error from finitely many replications. 
The condition $B_n/\log(n+1)\to\infty$ ensures this error vanishes uniformly over $n\geq m$ as $m\to\infty$.

For degenerate $U$-statistics of degree two, \citet{cai2026asymptotic} construct SAGE boundaries by combining
time-uniform Gaussian mean boundaries across Brownian coordinates, using eigenvalue estimates and coordinatewise significance-level allocation. 
Boot-AVT instead estimates a single quantile of the logarithmically penalized maximum of the full Gaussian chaos process through a quadratic multiplier bootstrap based on an empirical estimate of the second canonical function. 
This construction avoids explicit eigenvalue estimation and coordinatewise significance-level allocation. 
The simulations in Section~\ref{sec:simulations} show higher empirical power for the proposed procedures in the settings considered, while the real-data applications in Section~\ref{sec:realdata} show earlier rejection than most competing procedures.

We next study the power and stopping time of Boot-AVT under a fixed alternative. 
Let $\theta_1=\E\{h(X_1,\ldots,X_q)\}>\theta_0$ and $\Delta=\theta_1-\theta_0$. 
The following proposition bounds the probability of rejection by
sample size $N$ and establishes power one with a finite expected
stopping time.

\begin{proposition}\label{prop:fixed-alternative-consistency}
    Suppose Assumptions~\ref{ass:front-kernel} and~\ref{ass:front-multiplier} hold and $\theta=\theta_1>\theta_0$.
    Fix $\alpha\in(0,1)$ and $c>0$. 
    There exists a finite constant $C>0$ such
    that, for all integers $N\geq m\geq q$ satisfying
    $\Delta>(a_qc/N)\log(N/m)$,
    \begin{equation*}
        \pr(\tau_m^{(B_n)}\leq N)
        \geq1-
        \frac{C(N^{-p/2}+\alpha^{-1}N^{1-p})}
        {\{\Delta-(a_qc/N)\log(N/m)\}^{p}}.
    \end{equation*}
    Moreover, $\pr(\tau_m^{(B_n)}<\infty)=1$ and $\E(\tau_m^{(B_n)})<\infty$ for every fixed integer $m\geq q$.
\end{proposition}

\subsection{Sequential implementation}\label{sec:online-implementation}

Boot-AVT updates the $U$-statistic, the estimated second canonical function and the bootstrap maxima as observations accumulate. 
The kernel $h$, null value $\theta_0$, level $\alpha$, parameter $c$ and initial sample size $m\geq q$ are prespecified. 
For implementation, take $B_n$ to be nondecreasing. 
Bootstrap statistics are computed from sample size $q$ onward, and the rejection rule is applied from sample size $m$.

At sample size $q$, $U_{q,q}(h)=h(X_1,\ldots,X_q)$ and $\widehat h_q^{(2)}(X_i,X_j)=0$ for $1\leq i<j\leq q$. 
For each replication $b\leq B_q$, generate $\xi_1^{(b)},\ldots,\xi_q^{(b)}$ and initialize $\widetilde T_{q,2}^{*,b}=H_{c,q}^{*,b}=0$. 
When $m=q$, the initial critical value is zero and the rejection rule is applied immediately after initialization.

After observing $X_n$ for $n>q$, update the $U$-statistic by
\begin{equation}\label{eq:update-U}
    U_{n,q}(h) = \frac{n-q}{n}U_{n-1,q}(h)+\binom nq^{-1}\sum_{1\leq r_1<\cdots<r_{q-1}<n}h(X_{r_1},\ldots,X_{r_{q-1}},X_n).
\end{equation}
The averages $\widehat g_n^{(1)}$ and $\widehat g_n^{(2)}$ are
updated by separating terms involving $X_n$ from those based on
the first $n-1$ observations.
For $q>2$ and $1\leq i<j<n$,
\begin{equation}\label{eq:update-g2}
    \widehat g_n^{(2)}(X_i,X_j) = \frac{n-q}{n-2}\widehat g_{n-1}^{(2)}(X_i,X_j)+\binom{n-2}{q-2}^{-1}\sum_{\substack{1\leq r_1<\cdots<r_{q-3}<n \\ r_\ell\notin\{i,j\},\ 1\leq\ell\leq q-3}}h(X_i,X_j,X_{r_1},\ldots,X_{r_{q-3}},X_n).
\end{equation}
For $q=3$, the added term is $h(X_i,X_j,X_n)/(n-2)$.
For $q>2$ and $1\leq i<n$,
\begin{equation}\label{eq:update-g1}
    \widehat g_n^{(1)}(X_i) = \frac{n-q}{n-1}\widehat g_{n-1}^{(1)}(X_i)+\binom{n-1}{q-1}^{-1}\sum_{\substack{1\leq r_1<\cdots<r_{q-2}<n \\ r_\ell\ne i,\ 1\leq\ell\leq q-2}}
    h(X_i,X_{r_1},\ldots,X_{r_{q-2}},X_n).
\end{equation}
For $q=2$, these formulas reduce to $\widehat g_n^{(2)}(X_i,X_j)=h(X_i,X_j)$ and $\widehat g_n^{(1)}(X_i)=\{(n-2)\widehat g_{n-1}^{(1)}(X_i)+h(X_i,X_n)\}/(n-1)$.
The new values $\widehat g_n^{(1)}(X_n)$ and
$\widehat g_n^{(2)}(X_i,X_n)$, $1\leq i<n$, are computed directly
from \eqref{eq:empirical-conditional-kernel}.
Substitution into \eqref{eq:empirical-second-canonical} gives
$\widehat h_n^{(2)}(X_i,X_j)$ for all $1\leq i<j\leq n$.

For each existing replication $b\leq B_{n-1}$, generate
$\xi_n^{(b)}$ and update the current bootstrap statistic by
\begin{equation}\label{eq:sequential-bootstrap-current}
    \begin{aligned}
        \widetilde T_{n,2}^{*,b}
        ={}&\frac{n-2}{n-1}\widetilde T_{n-1,2}^{*,b}
        +\frac{2\xi_n^{(b)}}{n-1}
        \sum_{i=1}^{n-1}\xi_i^{(b)}\widehat h_n^{(2)}(X_i,X_n)\\
        &+\frac{2}{n-1}\sum_{1\leq i<j<n}\xi_i^{(b)}\xi_j^{(b)}
        \left\{\widehat h_n^{(2)}(X_i,X_j)
        -\widehat h_{n-1}^{(2)}(X_i,X_j)\right\}.
    \end{aligned}
\end{equation}
The second term adds pairs involving $X_n$, while the last term
accounts for changes in the estimated kernel at previously observed
pairs.
The logarithmic penalty changes by $-c\log\{n/(n-1)\}$ at every
earlier sample size, so
\begin{equation}\label{eq:sequential-bootstrap-recursion}
    H_{c,n}^{*,b}
    =\max\left\{H_{c,n-1}^{*,b}-c\log\{n/(n-1)\},
    \widetilde T_{n,2}^{*,b}\right\}.
\end{equation}
Thus each existing maximum is updated from its previous value and
the current bootstrap statistic, without recomputing
$\widetilde T_{k,2}^{*,b}$ for $k<n$.

When $B_n>B_{n-1}$, generate an independent multiplier sequence
$\xi_1^{(b)},\ldots,\xi_n^{(b)}$ for each new replication
$b=B_{n-1}+1,\ldots,B_n$. 
Starting from $\widetilde T_{q,2}^{*,b}=H_{c,q}^{*,b}=0$,
apply \eqref{eq:sequential-bootstrap-current} and
\eqref{eq:sequential-bootstrap-recursion} successively at sample
sizes $k=q+1,\ldots,n$, using $\widehat h_k^{(2)}(\cdot,\cdot)$ at sample size $k$.
These calculations use the stored estimates or the
observations from which they can be recomputed.

After updating all replications, compute
$\widehat\kappa_{1-\alpha,c}^{*,(B_n)}(n)$ as their empirical
$(1-\alpha)$-quantile and then the rejection rule 
$\phi_{m,n}^{(B_n)}(\theta_0) = \one\bigl\{T_n(\theta_0)-c\log(n/m)>\widehat\kappa_{1-\alpha,c}^{*,(B_n)}(n)\bigr\}$ with $T_n(\theta_0)=a_q^{-1}n\{U_{n,q}(h)-\theta_0\}$
when $n\geq m$.
Reject $H_0$ when $\phi_{m,n}^{(B_n)}(\theta_0)=1$.
Algorithm~\ref{alg:sequential-bootstrap-update} summarizes the sequential implementation of Boot-AVT.

\begin{algorithm}[htbp]
	\caption{Sequential implementation of Boot-AVT}
	\label{alg:sequential-bootstrap-update}
	\small
	\begin{algorithmic}[1]
		\Require Kernel $h$ of degree $q$, null value $\theta_0$, level
		$\alpha\in(0,1)$, parameter $c>0$, initial sample size $m\geq q$,
		and a nondecreasing sequence of positive integers $\{B_n\}_{n\geq q}$.
		\State Observe $X_1,\ldots,X_q$ and generate
		$\xi_1^{(b)},\ldots,\xi_q^{(b)}$ for $b=1,\ldots,B_q$.
		\State Initialize $U_{q,q}(h)$, $\widehat g_q^{(1)}(X_i)$ and
		$\widehat g_q^{(2)}(X_i,X_j)$ to $h(X_1,\ldots,X_q)$ for all
		$i$ and $i<j$. Set $\widehat h_q^{(2)}=0$ and
		$\widetilde T_{q,2}^{*,b}=H_{c,q}^{*,b}=0$ for $b=1,\ldots,B_q$.
		\If{$m=q$ and $T_q(\theta_0)>0$}
		\State Reject $H_0$ and stop.
		\EndIf
		\For{$n=q+1,q+2,\ldots$}
		\State Observe $X_n$ and update $U_{n,q}(h)$ by Equation~\eqref{eq:update-U}.
		\State Update $\widehat g_n^{(1)}$ and $\widehat g_n^{(2)}$ by Equations~\eqref{eq:update-g2} and \eqref{eq:update-g1} and compute $\widehat h_n^{(2)}$ by
		Equation~\eqref{eq:empirical-second-canonical}.
		\For{$b=1,\ldots,B_{n-1}$}
		\State Generate $\xi_n^{(b)}$ and update
		$\widetilde T_{n,2}^{*,b}$ and $H_{c,n}^{*,b}$ by
		Equations~\eqref{eq:sequential-bootstrap-current} and
		\eqref{eq:sequential-bootstrap-recursion}.
		\EndFor
		\For{$b=B_{n-1}+1,\ldots,B_n$}
		\State Generate an independent sequence
		$\xi_1^{(b)},\ldots,\xi_n^{(b)}$.
		\State Initialize $\widetilde T_{q,2}^{*,b}=H_{c,q}^{*,b}=0$.
		\State Apply Equations~\eqref{eq:sequential-bootstrap-current} and
		\eqref{eq:sequential-bootstrap-recursion} successively for
		$k=q+1,\ldots,n$, using $\widehat h_k^{(2)}$.
		\EndFor
		\State Compute $\widehat\kappa_{1-\alpha,c}^{*,(B_n)}(n)$.
		\If{$n\geq m$ and
			$T_n(\theta_0)-c\log(n/m)>\widehat\kappa_{1-\alpha,c}^{*,(B_n)}(n)$}
		\State Reject $H_0$ and stop.
		\EndIf
		\EndFor
	\end{algorithmic}
\end{algorithm}

\begin{remark}[Choice of $c$]\label{rem:choice-c}
    The parameter $c$ determines the trade-off between the initial critical value and the logarithmic term. 
    The critical value $\kappa_{1-\alpha,c}$ is nonincreasing in $c$, whereas $c\log(n/m)$ increases with $c$ for $n>m$.
    Consequently, increasing $c$ raises $\kappa_{1-\alpha,c}+c\log(n/m)$ at sufficiently late sample sizes.
    The choice of $c$ can therefore reflect the relative importance of early rejection and performance over longer monitoring periods. 
    One choice is $c=c_0s_T$, where $c_0>0$ is a small fixed coefficient and $s_T>0$ is a reference value for the standard deviation of the null limiting distribution of $T_n(\theta_0)$, 
    so that 
    \begin{equation*}
        \{T_n(\theta_0)-c\log(n/m)\}/{s_T}=\{T_n(\theta_0)\}/{s_T}-c_0\log(n/m).
    \end{equation*}
    Thus $c_0$ is the coefficient of the logarithmic term for $T_n(\theta_0)/s_T$. The rejection rule is invariant when $h-\theta_0$ and its reference scale $s_T$ are multiplied by the same positive constant.
    This reference value may be specified from prior knowledge,
    obtained by simulation under a reference null model, or estimated from an independent pilot sample. 
    The selected value of $c$ is held fixed throughout monitoring. 
\end{remark}

\section{Extensions to multiple kernels}
\label{sec:joint-finite-family}

Boot-AVT uses a single kernel chosen before monitoring. 
In practice, the choice of kernel, parameter or feature representation can affect power, and a suitable choice depends on the alternative. 
Boot-AVT-Multi combines a fixed finite collection of kernels in a joint test. 
Examples include Gaussian kernels with different bandwidths for
the MMD test in Example~\ref{ex:prelim-mmd} and different pairs of kernels for the HSIC test in Example~\ref{ex:prelim-hsic}. 
The collection may also include kernels applied to different fixed feature representations.

At a fixed sample size, \citet{SchrabEtAl2023} aggregate MMD tests over kernels and bandwidths, and \citet{ChatterjeeBhattacharya2025} combine MMD statistics using their joint asymptotic distribution and a multiplier bootstrap. 
In the sequential setting, the critical value must account for the maximum over both candidate statistics and all monitoring times. 
This requires accounting for the dependence among statistics computed from the same observations. 
Boot-AVT-Multi extends the time-inversion argument to the joint limiting process and estimates its critical values using common multipliers across kernels and sample sizes. 
The resulting joint test has asymptotic size $\alpha$, with the
critical value updated at each sample size.

\subsection{Joint limiting distribution}
\label{sec:joint-limit}

Let $J$ be fixed and finite. 
For $j=1,\ldots,J$, let $h_j:\mathcal X^{q_j}\to\mathbb R$
be a symmetric kernel of fixed degree $q_j\geq2$, and define
$\theta_j=\E\{h_j(X_1,\ldots,X_{q_j})\}$. 
We consider testing
\begin{equation*}
    H_0:\theta_j=\theta_{0,j},\quad j=1,\ldots,J,
    \qquad\text{versus}\qquad
    H_1:\theta_j>\theta_{0,j}\text{ for at least one }j.
\end{equation*}
At each sample size $n\geq q_m=\max_{1\leq j\leq J}q_j$, we estimate $\theta_j$ by $U_{n,q_j}(h_j)$, $j=1,\ldots,J$, all computed from $X_1,\ldots,X_n$. 
Let $h_j^{(1)}$ and $h_j^{(2)}$ denote the first two canonical
functions of $h_j$. 
We assume that, under $H_0$, each $U_{n,q_j}(h_j)$ is degenerate
of order two, that is, 
$h_j^{(1)}(X_1)=0$ almost surely and $0<\E[\{h_j^{(2)}(X_1,X_2)\}^2]<\infty$, $j=1,\ldots,J$.

Let $s_1,\ldots,s_J>0$ be fixed reference scales and
$a_j=q_j(q_j-1)/2$.
Define
\begin{equation}
    \notag
    T_{j,n}(\theta_j) = n\{U_{n,q_j}(h_j)-\theta_j\}/(a_js_j)
\end{equation}
for $n\geq q_j$. 
The scales place the candidate statistics on comparable numerical
ranges; they may be chosen using reference null models as in
Remark~\ref{rem:choice-c}.
For fixed $c>0$ and $n\geq m\geq q_m$, the joint statistic is
\begin{equation}
    \notag 
    A_{J,m,n} = \max_{1\leq j\leq J}T_{j,n}(\theta_{0,j})-c\log(n/m).
\end{equation}

For each $j$, write $\overline h_j^{(2)}=h_j^{(2)}/s_j$. 
Let $\lambda_{j\ell}$ and $\psi_{j\ell}$ denote its eigenvalues and orthonormal eigenfunctions, so that
$\overline h_j^{(2)}(x,y)=\sum_{\ell\geq1}\lambda_{j\ell}\psi_{j\ell}(x)\psi_{j\ell}(y)$
in $L^2(P\times P)$, with $\E\{\psi_{j\ell}(X_1)\}=0$ and
$\E\{\psi_{j\ell}(X_1)\psi_{jr}(X_1)\}=\one(\ell=r)$. 
The following assumptions extend the conditions in Section~\ref{sec:methodology} to the collection of kernels.

\begin{assumption}\label{ass:joint-spectrum}
    Under $H_0$, for every $j=1,\ldots,J$,  $0<\Lambda_{1,j}=\sum_{\ell\geq1}|\lambda_{j\ell}|<\infty$.
\end{assumption}

\begin{assumption}\label{ass:multi-kernel}
    $\E|h_j(X_1,\ldots,X_{q_j})|^p<\infty$ for every $j=1,\ldots,J$.
\end{assumption}

To describe the joint limit, let
$\{W_{j\ell}:1\leq j\leq J,\ \ell\geq1\}$ be a jointly Gaussian family of standard Brownian motions satisfying
\begin{equation}\label{eq:joint-brownian-covariance}
    \E\{W_{j\ell}(s)W_{j'\ell'}(t)\}
    =(s\wedge t)\E\{\psi_{j\ell}(X_1)\psi_{j'\ell'}(X_1)\},
    \qquad s,t\geq0.
\end{equation}
For each fixed $j$, $W_{j1},W_{j2},\ldots$ are independent, whereas the dependence across kernels is determined by the cross-kernel inner products in \eqref{eq:joint-brownian-covariance}. 
Thus each kernel has its own spectral representation, and these inner products enter the joint law in addition to the individual spectra. 
Define for $j=1,\ldots,J$,
\begin{equation*}
    Q_j(t)=\sum_{\ell\geq1}\lambda_{j\ell}\left\{W_{j\ell}^2(t)/t-1\right\},
    \qquad t>0.
\end{equation*}
Brownian scaling and time inversion hold jointly for this family. 
Consequently, for every $n\geq1$,
\begin{equation}
    \notag 
    \begin{aligned}
        \sup_{t\geq m}\max_{1\leq j\leq J}\{Q_j(t)-c\log(t/m)\}
        &\overset d=
        \sup_{0<s\leq1}\max_{1\leq j\leq J}\{Q_j(s)+c\log s\} \\
        &\overset d=
        \sup_{0<t\leq n}\max_{1\leq j\leq J}\{Q_j(t)+c\log(t/n)\}.
    \end{aligned}
\end{equation}
Similarly, let $H_{J,c}=\sup_{0<s\leq1}\max_{1\leq j\leq J}\{Q_j(s)+c\log s\}$ and $F_{J,c}(x)=\pr(H_{J,c}\leq x)$. 
Then its $(1-\alpha)$-quantile $\kappa_{J,1-\alpha,c}$ satisfies
\begin{equation*}
    \pr\left(
    \max_{1\leq j\leq J}Q_j(t)-c\log(t/m)\leq\kappa_{J,1-\alpha,c}
    \text{ for all }t\geq m
    \right)=1-\alpha.
\end{equation*}

For the $U$-statistics, define 
$J_{J,c,m}^{U}=\sup_{k\geq m}\max_{1\leq j\leq J}\{T_{j,k}(\theta_j)-c\log(k/m)\}$ for $m\geq q_m$ and 
$H_{J,c,n}^{U}=\max_{q_m\leq k\leq n}\max_{1\leq j\leq J}\{T_{j,k}(\theta_j)+c\log(k/n)\}$ for $n\geq q_m$. 
Under $H_0$, $\pr_0(J_{J,c,m}^{U}>\kappa)$ is the overall rejection probability when $A_{J,m,n}$ is compared with a fixed critical value $\kappa$, and $H_{J,c,n}^{U}$ is computed from the first $n$ observations by setting $\theta_j=\theta_{0,j}$.
Under Assumptions~\ref{ass:joint-spectrum} and~\ref{ass:multi-kernel},
\begin{equation}
    \notag
    H_{J,c,n}^{U}\xrightarrow{d}H_{J,c}\quad \text{as } n\to\infty,
    \qquad
    J_{J,c,m}^{U}\xrightarrow{d}H_{J,c}\quad \text{as } m\to\infty
\end{equation}
under $H_0$. 
This common limit allows the critical value for joint monitoring to be
estimated from the null distribution of $H_{J,c,n}^{U}$.

\subsection{Joint multiplier bootstrap}
\label{sec:joint-bootstrap}

To approximate the null distribution of $H_{J,c,n}^{U}$, we use
multipliers satisfying Assumption~\ref{ass:front-multiplier}.
Within each bootstrap replication, the same multiplier sequence
is used across sample sizes and candidate kernels.

Let $\widehat g_{j,k}^{(1)}$ and $\widehat g_{j,k}^{(2)}$ be the
averages in \eqref{eq:empirical-conditional-kernel} with $h$ and $q$
replaced by $h_j$ and $q_j$.
The corresponding estimate of $h_j^{(2)}$ is
$\widehat h_{j,k}^{(2)}(X_i,X_\ell)=\widehat g_{j,k}^{(2)}(X_i,X_\ell)-\widehat g_{j,k}^{(1)}(X_i)-\widehat g_{j,k}^{(1)}(X_\ell)+U_{k,q_j}(h_j)$.
For each bootstrap replication $b$, generate an independent copy
$\xi_1^{(b)},\xi_2^{(b)},\ldots$ of the multiplier sequence and define,
for $k\geq q_m$,
\begin{equation}
    \notag 
    \widetilde T_{j,k,2}^{*,b}
    =\frac{2}{(k-1)s_j}\sum_{1\leq i<\ell\leq k}
    \xi_i^{(b)}\xi_\ell^{(b)}\widehat h_{j,k}^{(2)}(X_i,X_\ell),
    \qquad j=1,\ldots,J.
\end{equation}
The bootstrap maximum is
\begin{equation}
    \notag
    H_{J,c,n}^{*,b}
    =\max_{q_m\leq k\leq n}\max_{1\leq j\leq J}
    \{\widetilde T_{j,k,2}^{*,b}+c\log(k/n)\}.
\end{equation}
The shared multipliers account for the dependence among the candidate statistics, as well as their dependence across sample sizes. 

Let $F_{J,n}^*(x)=\pr(H_{J,c,n}^{*,1}\leq x\mid X_1,\ldots,X_n)$. 
For a prespecified deterministic sequence of positive integers
$\{B_n\}_{n\geq q_m}$, define
$\widehat F_{J,c,n,B_n}^*(x)=B_n^{-1}\sum_{b=1}^{B_n}\one(H_{J,c,n}^{*,b}\leq x)$.
The estimated critical value is the empirical $(1-\alpha)$-quantile of the bootstrap maxima,
\begin{equation}
    \notag 
    \widehat\kappa_{J,1-\alpha,c}^{*,(B_n)}(n)
    =\inf\{x\in\mathbb R:\widehat F_{J,c,n,B_n}^*(x)\geq1-\alpha\}.
\end{equation}
The rejection rule and corresponding stopping time are
\begin{equation}
    \notag 
    \phi_{J,m,n}^{(B_n)}
    =\one\{A_{J,m,n}>\widehat\kappa_{J,1-\alpha,c}^{*,(B_n)}(n)\},
    \qquad
    \tau_{J,m}^{(B_n)}=\inf\{n\geq m:\phi_{J,m,n}^{(B_n)}=1\}.
\end{equation}
Taking $J=1$ and $s_1=1$ recovers the procedure in
Section~\ref{sec:methodology}. 
The following theorem gives bootstrap consistency and an explicit
bound on the overall type I error probability.

\begin{theorem}\label{thm:joint-asymptotic-size}
Under $H_0$, suppose Assumptions~\ref{ass:front-multiplier}, \ref{ass:joint-spectrum} and~\ref{ass:multi-kernel} hold. Fix $\alpha\in(0,1)$ and $c>0$, and suppose that $B_n/\log(n+1)\to\infty$. Let $\underline b_m=\inf_{n\geq m}\{B_n/\log(n+1)\}$. Then $\sup_{x\in\mathbb R}|F_{J,n}^*(x)-F_{J,c}(x)|\to0$ almost surely as $n\to\infty$, and $\sup_{n\geq m}|\widehat\kappa_{J,1-\alpha,c}^{*,(B_n)}(n)-\kappa_{J,1-\alpha,c}|\to0$ almost surely as $m\to\infty$. Moreover, there exists a finite constant $C_J>0$ such that, for all sufficiently large integers $m$,
\begin{align*}
&\left|\pr_0\left(A_{J,m,n}\leq\widehat\kappa_{J,1-\alpha,c}^{*,(B_n)}(n)\text{ for all }n\geq m\right)-(1-\alpha)\right|\\
\leq{}&C_J\{\log\log(m+1)\}^{1/2}\{\log(m+1)\}^{-p/2}+(2/\underline b_m)^{1/2}+2/(3m^3).
\end{align*}
Consequently,
\begin{equation*}
\begin{aligned}
&\pr_0\left(\exists~n\geq m:\phi_{J,m,n}^{(B_n)}=1\right)
=\pr_0(\tau_{J,m}^{(B_n)}<\infty)\\
=&\alpha+O\left(\{\log\log(m+1)\}^{1/2}\{\log(m+1)\}^{-p/2}+\underline b_m^{-1/2}\right).
\end{aligned}
\end{equation*}
\end{theorem}
Theorem~\ref{thm:joint-asymptotic-size} establishes asymptotic size $\alpha$ for joint monitoring. 
For fixed $J$, the bound has the same order in $m$ as in
Theorem~\ref{thm:front-final-bootstrap}; the constant $C_J$ now depends on the collection of kernels and their joint law. 
The remaining terms quantify the Monte Carlo error. 
Thus the procedure accounts for dependence across kernels through the joint bootstrap distribution.

\subsection{Power and sequential implementation}
\label{sec:joint-power-implementation}

Under a fixed alternative, let
$d_j=(\theta_j-\theta_{0,j})/(a_js_j)$ and
$d_*=\max_{1\leq j\leq J}d_j$.
The next theorem expresses the power bound in terms of the largest scaled discrepancy in the collection.

\begin{theorem}\label{thm:joint-fixed-alternative}
    Suppose Assumptions~\ref{ass:front-multiplier} and~\ref{ass:multi-kernel} hold under a fixed distribution satisfying $H_1$, so that $d_*>0$.
    Fix $\alpha\in(0,1)$ and $c>0$. Then there exists a finite constant $C_J>0$ such
    that, for all integers $N\geq m\geq q_m$ satisfying
    $d_*>(c/N)\log(N/m)$,
    \begin{equation*}
        \pr(\tau_{J,m}^{(B_n)}\leq N) \geq 1-\frac{C_J\{N^{-p/2}+\alpha^{-1}N^{1-p}\}}{\{d_*-(c/N)\log(N/m)\}^{p}}.
    \end{equation*}
    Moreover, for every fixed integer $m\geq q_m$,
    $\pr(\tau_{J,m}^{(B_n)}<\infty)=1$ and $\E(\tau_{J,m}^{(B_n)})<\infty$.
\end{theorem}

Thus a positive discrepancy for at least one candidate is sufficient for power one and a finite expected stopping time. 
The choice of the collection balances sensitivity to the alternative against the joint critical value: adding kernels may increase both $d_*$ and $\kappa_{J,1-\alpha,c}$. 
A small grid of bandwidths around a reference value gives one practical choice, as in the simulations in Section~\ref{sec:simulations}. 
The kernels, reference scales $s_j$ and parameter $c$ are fixed before monitoring; they may be chosen using prior knowledge or an independent pilot sample.

For sequential implementation, take $B_n$ to be nondecreasing. 
At the first common sample size $q_m$, compute the candidate
bootstrap statistics and initialize
$H_{J,c,q_m}^{*,b}=\max_{1\leq j\leq J}\widetilde T_{j,q_m,2}^{*,b}$.
After observing $X_n$ for $n>q_m$, append the same multiplier
$\xi_n^{(b)}$ to all candidates in replication $b$ and update their bootstrap statistics using \eqref{eq:sequential-bootstrap-current} applied to $h_j/s_j$.
The joint maximum then satisfies
\begin{equation}
    \notag 
    H_{J,c,n}^{*,b}
    =\max\left\{
    H_{J,c,n-1}^{*,b}-c\log\{n/(n-1)\},
    \max_{1\leq j\leq J}\widetilde T_{j,n,2}^{*,b}
    \right\}.
\end{equation}
When $B_n>B_{n-1}$, each new replication uses a multiplier sequence of length $n$ shared by all kernels; its historical statistics are computed for $q_m\leq k\leq n$ using the corresponding $\widehat h_{j,k}^{(2)}$.
After updating the replications, compute
$\widehat\kappa_{J,1-\alpha,c}^{*,(B_n)}(n)$ and $\phi_{J,m,n}^{(B_n)}$ whenever $n\geq m$. 

\section{Simulation results}\label{sec:simulations}

We examine the finite-sample performance of \textit{Boot-AVT} in goodness-of-fit, two-sample, and independence testing, and evaluate \textit{Boot-AVT-Multi} in the latter two problems. The experiments assess cumulative type I error, power as observations accumulate, and power at a fixed monitoring horizon over a range of signal strengths.

All tests use nominal level $\alpha=0.1$ and monitoring start times $m_0\in\{200,400,600\}$. 
The maximum sample size is $N=2000$. 
We use $R=2000$ independent replications in each setting. 
For each method, let $\tau^{(r)}(m_0)$ denote the first rejection time in replication $r$. We report the cumulative rejection probability
$\widehat\pi(n;m_0)=R^{-1}\sum_{r=1}^{R}\one\{\tau^{(r)}(m_0)\leq n\}$
for $m_0\leq n\leq N$. 
Under $H_0$, this estimates the probability of at least one rejection between $m_0$ and $n$; under $H_1$, it estimates power by time $n$. 

Boot-AVT and Boot-AVT-Multi use standard normal multipliers and $B_n=B=500$ bootstrap replications at each sample size. 
We set $c=0.01s_T$ for Boot-AVT, where $s_T$ is a reference standard deviation under the null model, and $c=0.01$ for Boot-AVT-Multi after scaling each candidate statistic by its reference standard deviation $s_j$. 
The MMD and HSIC kernels are Gaussian,
$K_\sigma(u,v)=\exp\{-\|u-v\|^2/(2\sigma^2)\}$.
Single-kernel procedures use $\sigma^2=1$. 
Boot-AVT-Multi uses $J=5$ kernels with $\sigma_j^2\in\{1/4,1/2,1,2,4\}$; for HSIC, each candidate uses the same bandwidth for both marginal kernels. 

\subsection{Goodness-of-fit tests}\label{sec:gof}

We observe $X_i$ independently from $\mathcal N(\delta,1)$ and test $H_0:F=\Phi$ against $H_1:F\neq\Phi$, where $F$ is the observation distribution function and $\Phi$ is the standard normal distribution function. The null corresponds to $\delta=0$, and positive values of $\delta$ give location alternatives. Boot-AVT uses the unbiased Cram\'er--von Mises statistic in Example~\ref{ex:prelim-gof}, with kernel $h(x,y)=\{\Phi(x)^2+\Phi(y)^2\}/2-\max\{\Phi(x),\Phi(y)\}+1/3$ and statistic $T_n(0)=nU_{n,2}(h)$.

We compare Boot-AVT with SAGE-GM and SAGE-LIL \citep{cai2026asymptotic}, the kernel betting test (Bet-MMD) \citep{shekhar2024nonparametric}, the MMD confidence sequence method (MR-MMD) \citep{manole2023martingale}, the sequential Kolmogorov--Smirnov test (HR-KS) \citep{howard2022sequential}, the kernelized Stein discrepancy betting test (MTR-KSD) \citep{martineztaboada2025sequential}, and the fixed-sample Cram\'er--von Mises test (CvM) \citep{csorgo1996exact}. 

\begin{figure}[!htb]
\centering
  \includegraphics[width=0.9\linewidth]{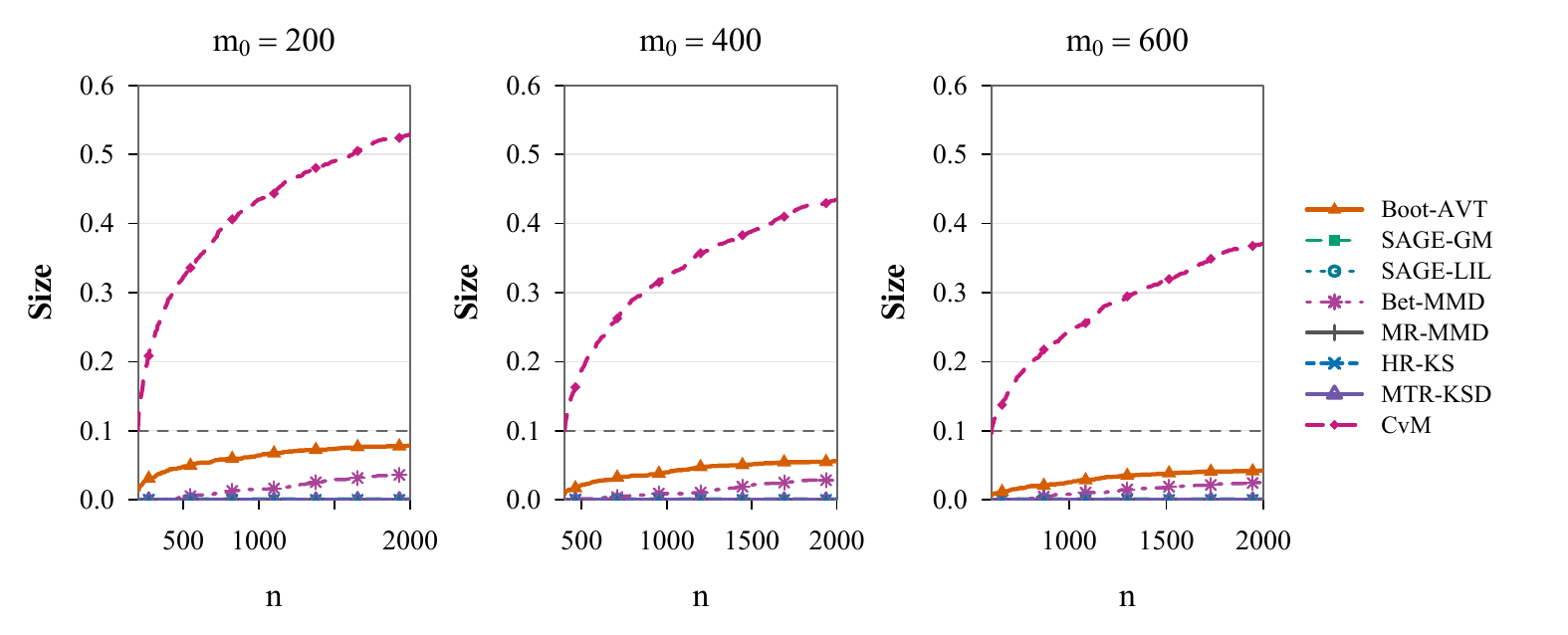}
\caption{Empirical size of the goodness-of-fit tests for different monitoring start times $m_0$.}
\label{fig:gof-anytime-size-all-methods}
\end{figure}

Figure~\ref{fig:gof-anytime-size-all-methods} reports cumulative rejection probabilities under $\delta=0$. Boot-AVT remains below the nominal level for all three monitoring start times, while the other sequential procedures have lower rejection probabilities. Repeated application of CvM produces substantial size inflation.

\begin{figure}[!htb]
\centering
  \includegraphics[width=0.9\linewidth]{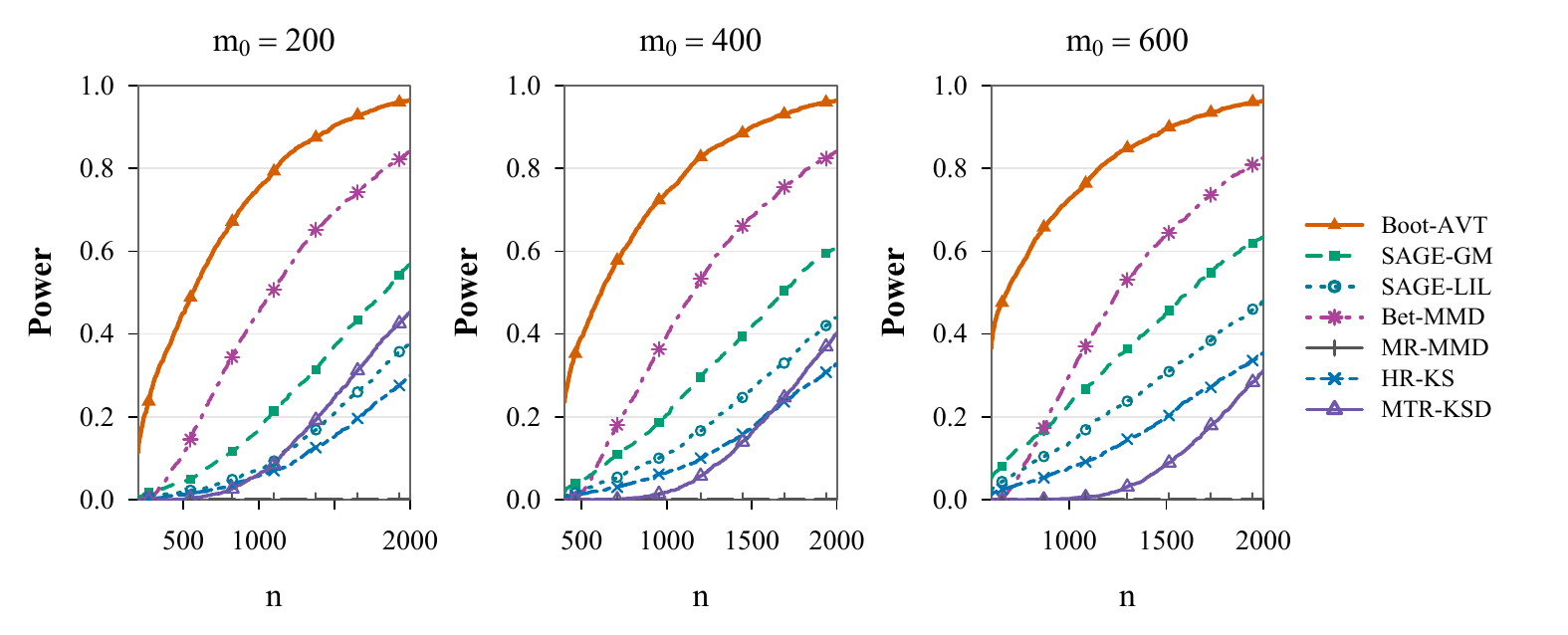}
\caption{Empirical power of the goodness-of-fit tests under the fixed Gaussian location shift $\delta=0.1$.}
\label{fig:gof-fixed-signal-power-all-methods}
\end{figure}

Under $\delta=0.1$, Figure~\ref{fig:gof-fixed-signal-power-all-methods} shows that Boot-AVT has higher cumulative power than the competing sequential procedures across the three monitoring start times. 
Bet-MMD is the strongest competitor under this alternative, whereas MR-MMD records no rejections over the monitoring period. 

\begin{figure}[!htb]
	\centering
	\includegraphics[width=0.9\linewidth]{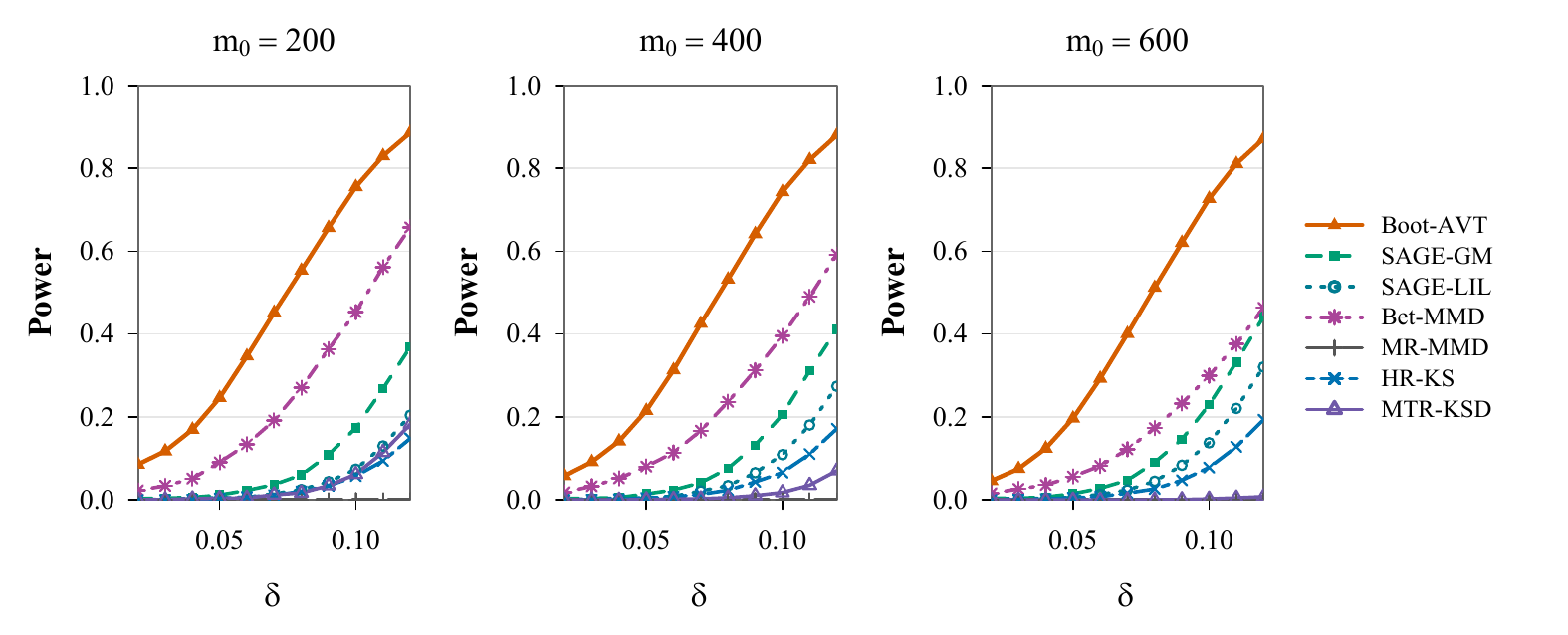}
	\caption{Empirical power of the goodness-of-fit tests at sample size $n=1000$.}
	\label{fig:gof-fixed-horizon-power-all-methods}
\end{figure}

Figure~\ref{fig:gof-fixed-horizon-power-all-methods} compares power by fixing the sample size at $n=1000$ and varying the signal strength $\delta$. Boot-AVT retains its power advantage over the competing sequential procedures, with similar performance across the three monitoring start times.

\subsection{Two-sample tests}
\label{sec:degenerate-mmd}

At time $i$, we independently observe $X_i=(Y_i,Z_i)$, where
$Y_i\sim P_Y=\mathcal N(0,1)$ and $Z_i\sim P_Z=\mathcal N(\delta,1)$,
and the two sequences are independent. 
We test $H_0:P_Y=P_Z$ against $H_1:P_Y\neq P_Z$; the null corresponds to $\delta=0$. 
Here and in the independence experiment, $n$ counts observed pairs. Boot-AVT uses $T_n(0)=nU_{n,2}(h_\sigma)$, where the paired MMD kernel is
$h_\sigma(X_i,X_j)=K_\sigma(Y_i,Y_j)+K_\sigma(Z_i,Z_j)-K_\sigma(Y_i,Z_j)-K_\sigma(Y_j,Z_i)$.
Boot-AVT-Multi jointly monitors the five bandwidths specified above.

We compare Boot-AVT and Boot-AVT-Multi with SAGE-GM and SAGE-LIL \citep{cai2026asymptotic}, Bet-MMD \citep{shekhar2024nonparametric}, the prediction-based betting test (Seq-C-2ST) \citep{podkopaev2023predictive}, MR-MMD \citep{manole2023martingale}, and the fixed-sample MMD permutation test (MMDp) \citep{gretton2012kernel}. 

\begin{figure}[!htb]
\centering
  \includegraphics[width=0.9\linewidth]{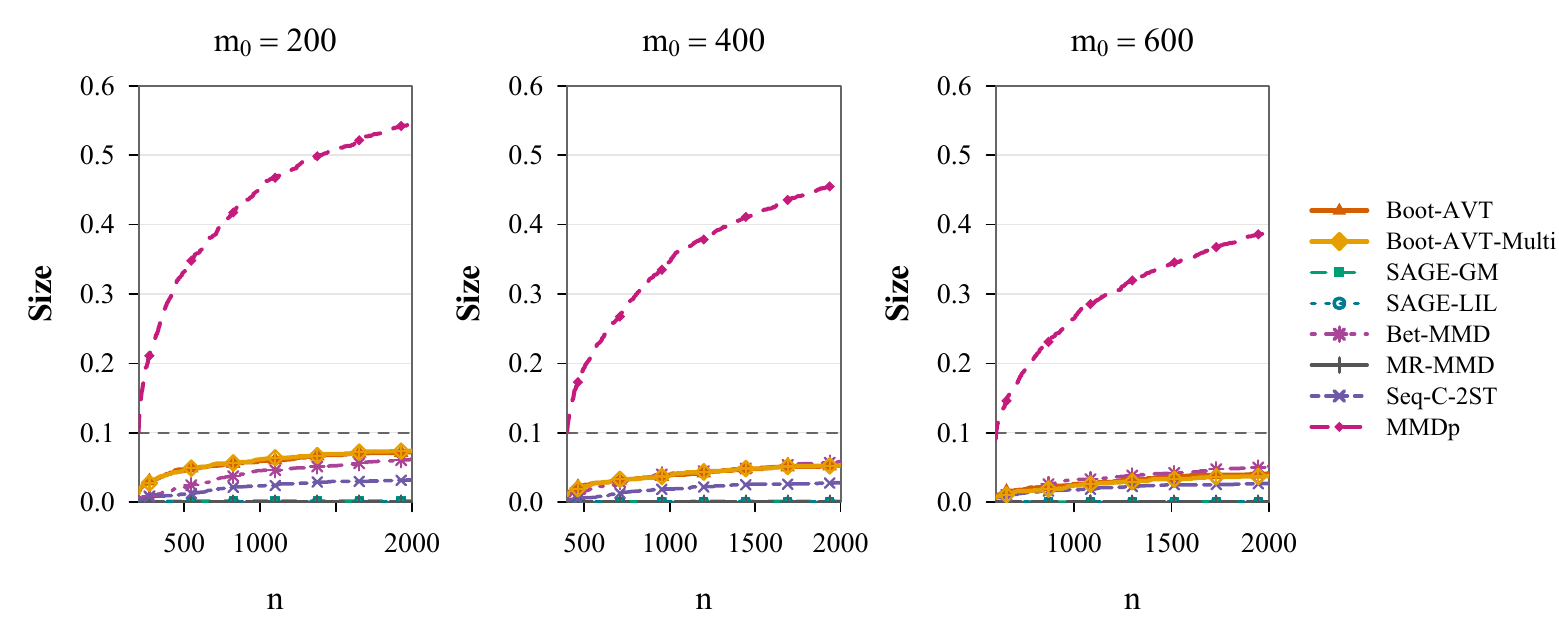}
\caption{Empirical size of the two-sample tests for different monitoring start times $m_0$.}
\label{fig:mmd-anytime-size-all-methods}
\end{figure}

Figure~\ref{fig:mmd-anytime-size-all-methods} shows empirical size below the nominal level for Boot-AVT and Boot-AVT-Multi at all three monitoring start times. 
The other sequential procedures are more conservative, while repeated application of MMDp produces size inflation.

\begin{figure}[!htb]
\centering
  \includegraphics[width=0.9\linewidth]{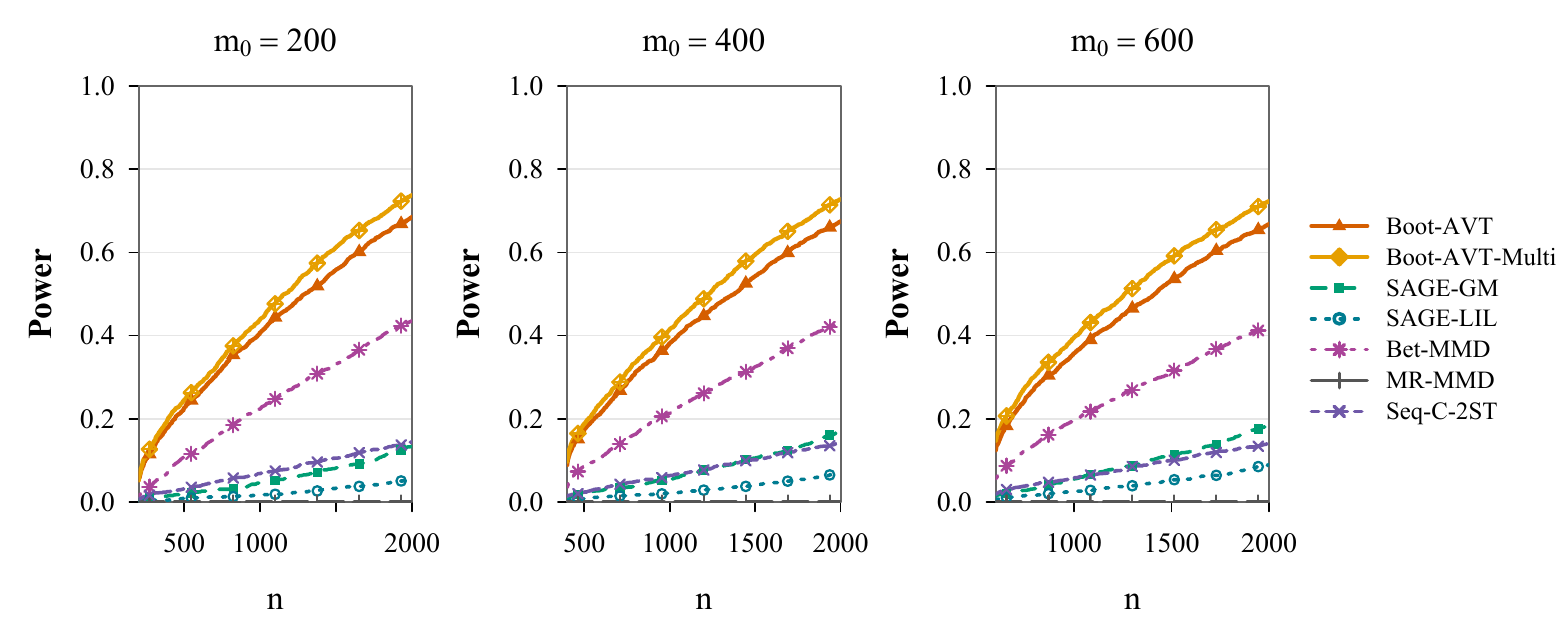}
\caption{Empirical power of the two-sample tests under the fixed Gaussian location shift $\delta=0.1$.}
\label{fig:mmd-fixed-signal-power-all-methods}
\end{figure}

Under $\delta=0.1$, Figure~\ref{fig:mmd-fixed-signal-power-all-methods} shows that Boot-AVT and Boot-AVT-Multi have higher cumulative power than the competing sequential procedures. 
Boot-AVT-Multi provides a modest improvement over Boot-AVT across the three monitoring start times. 
Bet-MMD is the strongest competitor, with power between $0.4$ and $0.45$ by $n=2000$. SAGE-GM and Seq-C-2ST have terminal power of approximately $0.1$--$0.2$, SAGE-LIL remains below $0.1$, and MR-MMD records no rejections. 

\begin{figure}[!htb]
	\centering
	\includegraphics[width=0.9\linewidth]{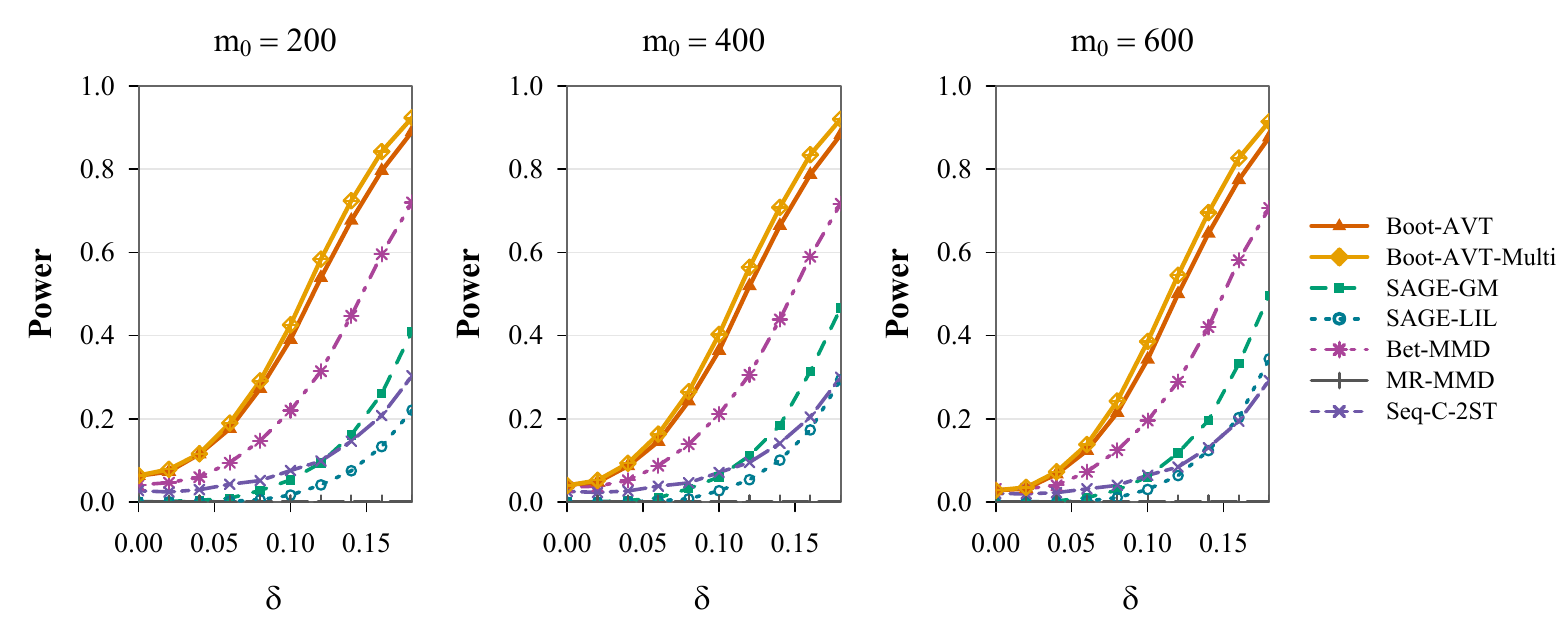}
	\caption{Empirical power of the two-sample tests at sample size $n=1000$.}
	\label{fig:mmd-fixed-horizon-power-all-methods}
\end{figure}

The comparison across location shifts leads to the same ordering, as shown in Figure~\ref{fig:mmd-fixed-horizon-power-all-methods}.  Boot-AVT and Boot-AVT-Multi have similar power by $n=1000$ and exceed the competing sequential procedures over the examined alternatives.

\subsection{Independence tests}
\label{subsec:independence-tests}

We observe independent pairs $X_i=(Y_i,Z_i)$ generated by
$Y_i,\varepsilon_i \sim \mathcal N(0,1)$ and 
$Z_i=\rho Y_i+\sqrt{1-\rho^2}\,\varepsilon_i$, $0\leq\rho<1$,
where the $Y_i$ and $\varepsilon_i$ sequences are independent. We test $H_0:P_{YZ}=P_Y\otimes P_Z$ against $H_1:P_{YZ}\neq P_Y\otimes P_Z$. The null corresponds to $\rho=0$; increasing $\rho$ introduces dependence while preserving standard normal marginal distributions.
Boot-AVT uses the unbiased HSIC statistic in Example~\ref{ex:prelim-hsic} \citep{SongEtAl2012}.
Both marginal kernels are Gaussian with $\sigma^2=1$. Boot-AVT-Multi uses the five bandwidths and candidate scales specified above.

We compare Boot-AVT and Boot-AVT-Multi with the HSIC betting tests SKIT-ONS and SKIT-aGRAPA \citep{podkopaev2023skit}, the rank-based sequential independence test (Sinkhorn) \citep{henzi2024rank}, the classification-based betting test (Seq-C-IT) \citep{podkopaev2023predictive}, the sequential independence test (SR-KS) \citep{shekhar2024nonparametric}, and the fixed-sample HSIC permutation test (HSICp) \citep{gretton2008kernel}. 

\begin{figure}[!htb]
\centering
  \includegraphics[width=0.9\linewidth]{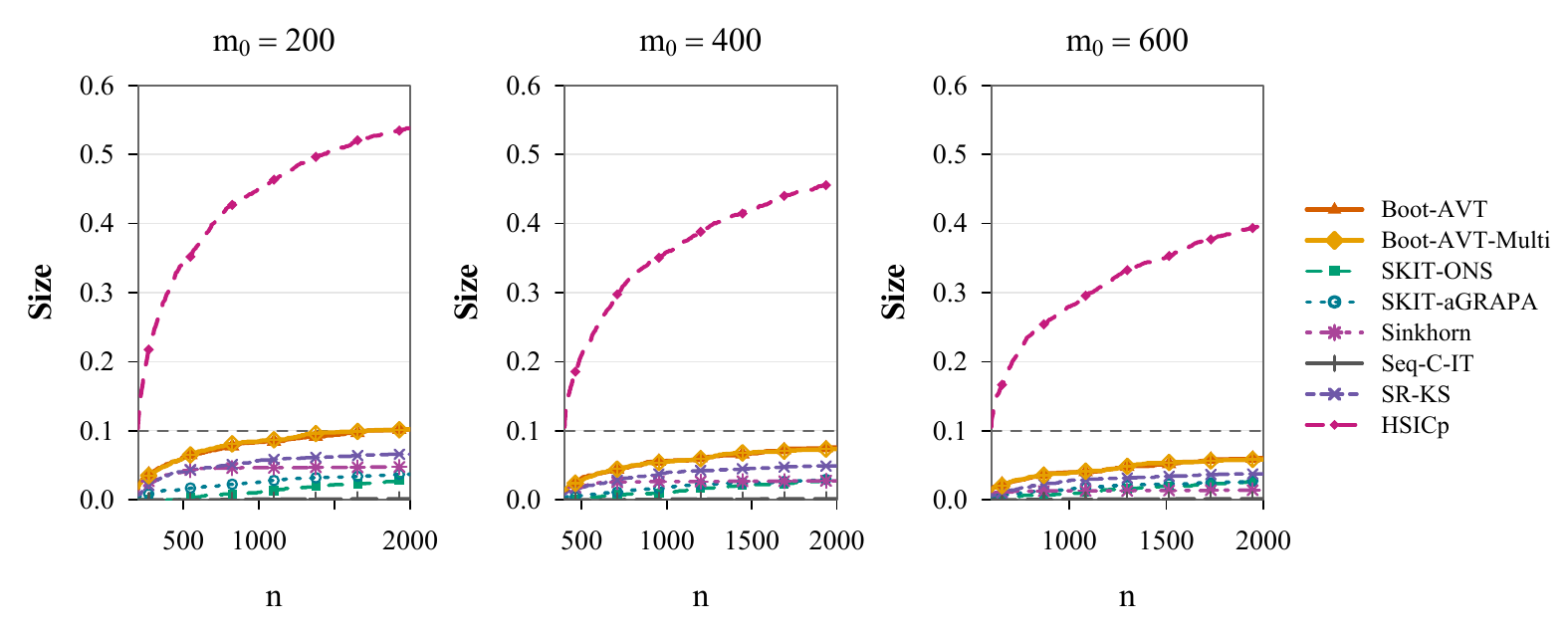}
\caption{Empirical size of the independence tests for different monitoring start times $m_0$.}
\label{fig:independence-size}
\end{figure}

In Figure~\ref{fig:independence-size}, Boot-AVT and Boot-AVT-Multi have cumulative rejection probabilities close to the nominal level for $m_0=200$ and below it for the later monitoring starts. 
The competing sequential procedures have smaller rejection probabilities. 
Repeated application of HSICp again produces substantial size inflation.

\begin{figure}[!htb]
\centering
  \includegraphics[width=0.9\linewidth]{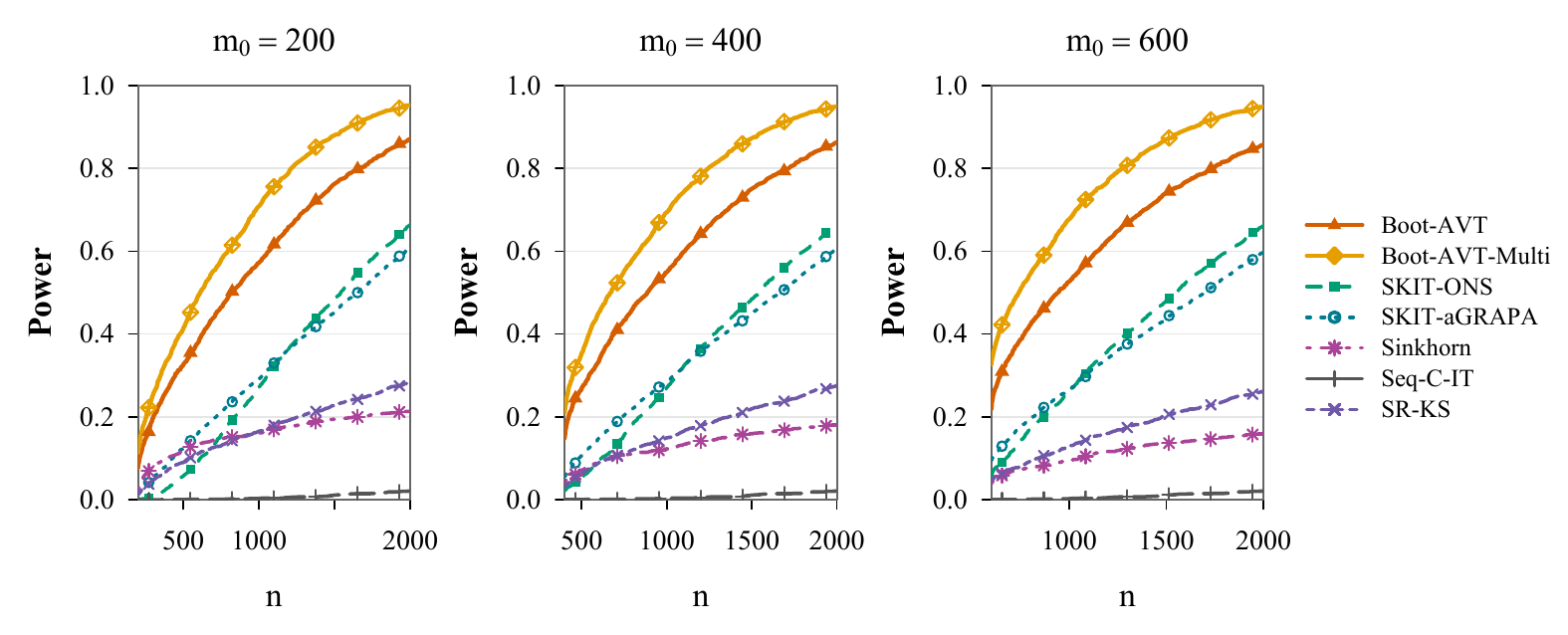}
\caption{Empirical power of the independence tests under the fixed Gaussian correlation $\rho=0.1$.}
\label{fig:independence-power-fixed-signal}
\end{figure}

Under $\rho=0.1$, Figure~\ref{fig:independence-power-fixed-signal} shows that Boot-AVT-Multi achieves the highest cumulative power, followed by Boot-AVT, across the three monitoring start times. 
SKIT-ONS and SKIT-aGRAPA are the strongest competitors. Sinkhorn, Seq-C-IT, and SR-KS have lower power under this alternative. 

\begin{figure}[!htb]
	\centering
	\includegraphics[width=0.9\linewidth]{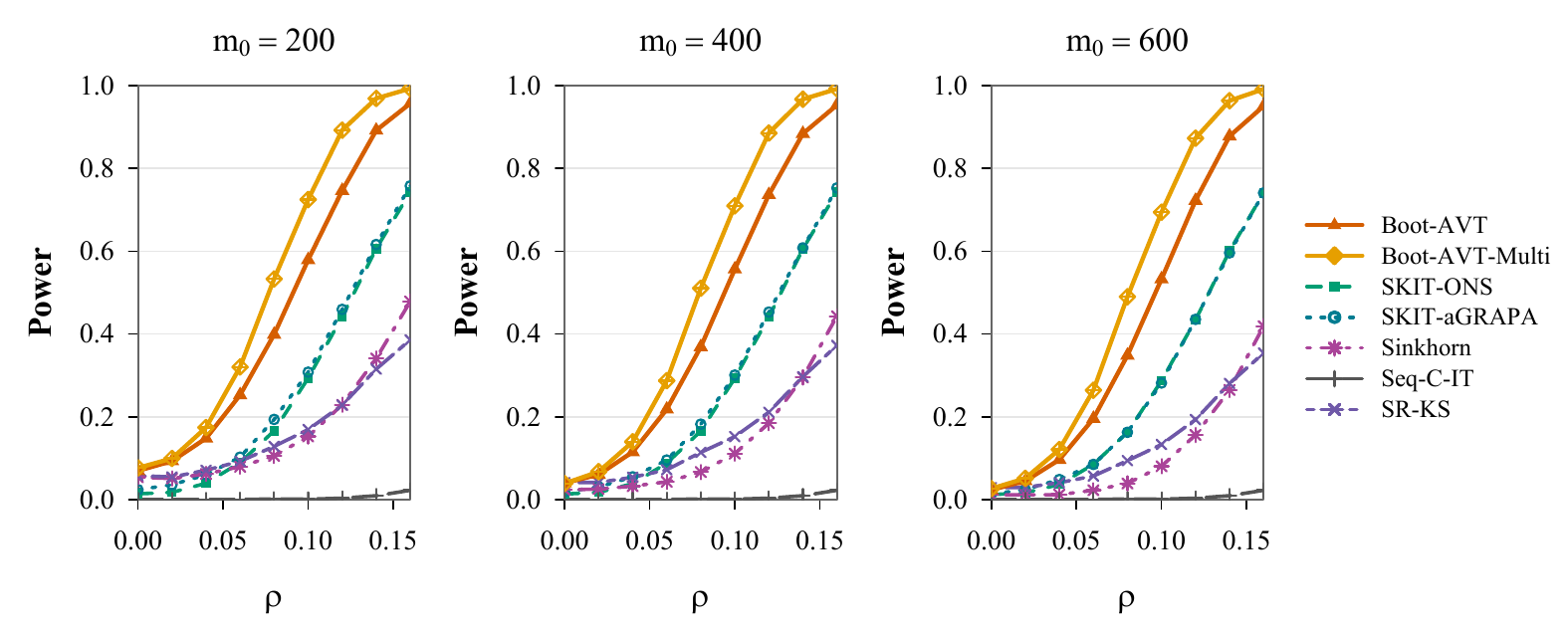}
	\caption{Empirical power of the independence tests at sample size $n=1000$.}
	\label{fig:independence-power-by-signal}
\end{figure}

Figure~\ref{fig:independence-power-by-signal} varies $\rho\in\{0,0.02,\ldots,0.16\}$ and reports power by $n=1000$. Boot-AVT-Multi improves on Boot-AVT over the intermediate signal strengths, and both procedures have higher power than the competing sequential tests.  The gain from combining bandwidths is more pronounced in this independence experiment than in the Gaussian two-sample location experiment.

\section{Real data applications}
\label{sec:realdata}

We apply Boot-AVT and Boot-AVT-Multi to a paired two-sample problem and an independence problem. 
Both procedures use nominal level $\alpha=0.1$, standard normal multipliers, and $B=500$ bootstrap replications. 
We set $c=0.01s_T$ for Boot-AVT and $c=0.01$ for Boot-AVT-Multi. 
In each application, the first 50 observations determine the Gaussian bandwidths and are excluded from subsequent testing. 
Boot-AVT-Multi uses $J=5$ kernels whose bandwidth variances equal the selected reference variance multiplied by $1/4$, $1/2$, $1$, $2$, and $4$. 
For HSIC, the same multiplier is applied to both marginal bandwidth variances. 
The competing procedures follow the settings in Section~\ref{sec:simulations}. 
All reported rejection times count observations in the testing sample.

\subsection{Right- versus left-ear audiometry}
\label{subsec:realdata-audiometry-m0200}

We analyze audiometry data from the 2015--2016 National Health and Nutrition Examination Survey (NHANES) \citep{NCHS2019Audiometry}. The analysis includes adults aged 20--69 years with pure-tone air-conduction thresholds recorded in both ears at 0.5, 1, 2, 3, 4, 6, and 8 kHz. 
For participant $i$, let $Y_i$ and $Z_i$ be the seven-dimensional vectors of right- and left-ear thresholds, respectively. 
We use the first of the repeated 1-kHz measurements. 
Removing observations with missing values leaves 4004 pairs. 
The paired MMD statistic tests equality of the ear-specific marginal distributions while retaining within-participant dependence:
\begin{equation*}
    H_0:P_Y=P_Z
    \qquad\text{versus}\qquad
    H_1:P_Y\neq P_Z.
\end{equation*}

We select a common reference bandwidth for Boot-AVT, Boot-AVT-Multi, SAGE-GM, SAGE-LIL, Bet-MMD, and MR-MMD using the first 50 pairs. 
The pooled median heuristic \citep{gretton2012kernel} takes the median of the positive pairwise Euclidean distances among the 50 right-ear and 50 left-ear vectors, giving $\widehat\sigma_{\mathrm{med}}=33.166$ dB HL. 
The remaining $N=3954$ pairs form the testing sample, and the test statistics follow Section~\ref{sec:degenerate-mmd}.

Table~\ref{tab:realdata-audiometry-m0200} shows that Boot-AVT and Boot-AVT-Multi reject at $n=231$ and $n=265$, respectively, when monitoring starts at $m_0=200$. 
Both reject immediately when monitoring starts at $m_0=400$ or $600$. 
Bet-MMD rejects at $n=741$, whereas SAGE-GM, SAGE-LIL, and Seq-C-2ST require more than 1600 pairs; MR-MMD does not reject by $N=3954$. 
The evidence of different ear-specific distributions is consistent with the hearing asymmetry reported by \citet{SuenEtAl2021} using earlier NHANES data.

\begin{table}[!htb]
\centering
\caption{First rejection times for the NHANES comparison of right- and left-ear threshold distributions. ``NR'' denotes no rejection by $N=3954$.}
\label{tab:realdata-audiometry-m0200}
\small
\setlength{\tabcolsep}{3pt}
\begin{tabular}{@{}rccccccc@{}}
\toprule
$m_0$ & Boot-AVT & Boot-AVT-Multi & SAGE-GM & SAGE-LIL & Bet-MMD & Seq-C-2ST & MR-MMD \\
\midrule
200 & 231 & 265 & 1714 & 1837 & 741 & 1865 & NR \\
400 & 400 & 400 & 1701 & 1806 & 741 & 1865 & NR \\
600 & 600 & 600 & 1674 & 1785 & 741 & 1865 & NR \\
\bottomrule
\end{tabular}
\end{table}

\subsection{Age at diagnosis and prognostic severity in METABRIC}
\label{subsec:realdata-metabric-independence}

We next examine clinical data from the Molecular Taxonomy of Breast Cancer International Consortium (METABRIC), which contains genomic and clinical information on primary breast cancer \citep{curtis2012metabric,pereira2016metabric}. 
We investigate the association between age at diagnosis and the Nottingham Prognostic Index (NPI), a measure of breast cancer prognosis \citep{haybittle1982npi}. We use the patient-level clinical data available through cBioPortal \citep{gao2013cbioportal}. 
Let $Y_i$ denote age at diagnosis and $Z_i$ the supplied NPI value for patient $i$. 
The clinical file contains 2509 distinct patients; excluding observations with missing age or NPI values leaves 2287 patients. We test
\begin{equation*}
    H_0:P_{YZ}=P_Y\otimes P_Z
    \qquad\text{versus}\qquad
    H_1:P_{YZ}\neq P_Y\otimes P_Z.
\end{equation*}

The first 50 patients determine the reference Gaussian bandwidths for Boot-AVT, Boot-AVT-Multi, SKIT-ONS, and SKIT-aGRAPA. 
Applying the median heuristic separately to the positive pairwise distances for age and NPI gives bandwidths of $15.760$ years and $1.008$ NPI points, respectively. 
All methods are applied to the remaining $N=2237$ patients using the test statistics in Section~\ref{subsec:independence-tests}.

Table~\ref{tab:realdata-metabric-independence} reports evidence of dependence between age at diagnosis and NPI. 
For $m_0=200$ and $400$, Boot-AVT-Multi rejects earlier than Boot-AVT; the ordering reverses at $m_0=600$. 
Sinkhorn gives the earliest rejection at $m_0=200$. 
The earliest rejections at the other two monitoring start times are obtained by Boot-AVT-Multi for $m_0=400$ and by Boot-AVT for $m_0=600$. 
Both bootstrap procedures reject before SKIT-ONS and SKIT-aGRAPA for all three monitoring start times. 
Seq-C-IT and SR-KS do not reject by $N=2237$.

\begin{table}[!htb]
\centering
\caption{First rejection times for the METABRIC independence test of age at diagnosis and NPI. ``NR'' denotes no rejection by $N=2237$.}
\label{tab:realdata-metabric-independence}
\small
\setlength{\tabcolsep}{3pt}
\begin{tabular}{@{}rccccccc@{}}
\toprule
$m_0$ & Boot-AVT & Boot-AVT-Multi & SKIT-ONS & SKIT-aGRAPA & Seq-C-IT & SR-KS & Sinkhorn \\
\midrule
200 & 375 & 313 & 1664 & 1714 & NR & NR & 246 \\
400 & 501 & 417 & 1664 & 1714 & NR & NR & 481 \\
600 & 601 & 714 & 1664 & 1714 & NR & NR & 1499 \\
\bottomrule
\end{tabular}
\end{table}

\section{Concluding remarks}

This paper develops Boot-AVT for asymptotically anytime-valid testing with degenerate $U$-statistics and Boot-AVT-Multi for joint testing with a fixed finite collection of kernels. 
Brownian scaling and time inversion relate the distribution of the logarithmically penalized maximum of the limiting Gaussian chaos process over an unbounded monitoring horizon to that of a maximum on a bounded interval. 
A quadratic multiplier bootstrap based on the empirical second canonical function estimates the critical values without explicit eigenvalue estimation, using common multipliers to account for dependence across sample sizes and kernels. 
We establish uniform consistency of the critical values and asymptotic size $\alpha$ as the monitoring start tends to infinity, with explicit bounds separating statistical approximation and Monte Carlo errors. 
Under fixed alternatives, we derive finite-horizon power bounds and establish power one with finite expected stopping times for every fixed monitoring start. 
Simulations and data applications illustrate the procedures in goodness-of-fit, two-sample and independence testing. 
Future directions include extending the calibration to $U$-statistics with higher-order degeneracy and developing procedures for simultaneous monitoring of multiple data streams with false discovery rate control.

\bibliographystyle{chicago}
\bibliography{ref}

\end{document}